\documentclass[11pt,a4paper]{amsart}
\usepackage[T1]{fontenc}
\usepackage{lmodern}
\usepackage[margin=26mm]{geometry}
\usepackage{amsmath,amssymb,amsthm,mathtools}
\usepackage{microtype,needspace}
\usepackage{enumitem,booktabs,array}
\usepackage[colorlinks=true,linkcolor=blue,citecolor=blue,urlcolor=blue]{hyperref}
\setlist[enumerate]{label=\textup{(\roman*)},leftmargin=*,itemsep=2pt}
\newtheorem{theorem}{Theorem}[section]
\newtheorem{proposition}[theorem]{Proposition}
\newtheorem{lemma}[theorem]{Lemma}
\newtheorem{corollary}[theorem]{Corollary}
\theoremstyle{definition}
\newtheorem{definition}[theorem]{Definition}
\newtheorem{example}[theorem]{Example}
\theoremstyle{remark}
\newtheorem{remark}[theorem]{Remark}
\DeclareMathOperator{\Min}{Min}
\DeclareMathOperator{\Spec}{Spec}
\DeclareMathOperator{\Frac}{Frac}

\newcommand{\MP}[2]{\mathsf{MP}_{#1,#2}}
\newcommand{\Z}{\mathbb Z}

\newcommand{\Q}{\mathbb Q}
\newcommand{\C}{\mathbb C}

\numberwithin{equation}{section}
\allowdisplaybreaks[2]
\title[Spectral and arithmetic classifications of mpdf ideals]
{Spectral and arithmetic classifications of mixed-power difference factor absorbing ideals}

\date{}
\subjclass[2020]{Primary 13A15; Secondary 13B25, 13B30, 13F05, 13F30}
\keywords{Mixed-power difference factor absorbing ideal, minimal prime, root of unity,
nilpotent extension, idealization, discrete valuation ring, Dedekind domain}

\begin{document}
\raggedbottom
\author{Omar H. Taha}
\address{Department of Mathematics, College of Computer Science and Mathematics, Tikrit University, Tikrit, Iraq.}
\email{omar.h.tahamm2314@st.tu.edu.iq}

\author{Akram Salim}
\address{Department of Mathematics, College of Computer Science and Mathematics, Tikrit University, Tikrit, Iraq.}
\email{akr\_tel@tu.edu.iq}

\author{Hwankoo Kim}
\address{Division of Computer Engineering, Hoseo University, Republic of Korea.}
\email{hkkim@hoseo.edu}

\begin{abstract}
For integers $n\ge2$ and $k\ge1$, a proper ideal is $(kn,n)$-mixed-power
difference factor absorbing if membership of $a^{kn}-b^n$ forces membership
of $a^k-b$ or of the complementary geometric-sum factor. Extending the
known cubic spectral criterion, we classify radical ideals for arbitrary
exponents in terms of roots of unity in the fraction fields of their
minimal-prime quotients, without finiteness or comaximality assumptions.
The condition is independent of $k$ for radical ideals. We also extend
the finite comaximal-prime criterion for quartic pair absorption to
arbitrary radical ideals and identify its distinct three-component
obstruction. More generally,
requiring one of the two factors to belong to the radical is equivalent
to the radical itself satisfying the mixed-power condition for every
multiplier. This yields radical descent and extends the cubic primary
descent theorem. We also prove equivalence under polynomial and formal
power series extension for radical ideals. These
results yield criteria for pullbacks of domains and, for geometrically
reduced algebras, a characterization of geometric integrality by stability
under field extension when the characteristic does not divide $n$.
For $n\ge3$ in the same characteristics, square-zero extensions of
geometrically reduced algebras have this stability precisely when the
quotient is geometrically integral and the kernel is torsion-free. In characteristic
$p$ dividing $n$, the condition depends only on the $p$-part of $n$ and is
equivalent to a power-closure condition. We also classify square-zero
extensions of domains, including nonsplit extensions and arbitrary
idealization ideals with prime projection. For homogeneous idealizations
over $P$-primary bases with $2\notin P$, a related quartic criterion
requires no restriction on $3$.
Finally, we determine the exact
valuation criterion when $n$ is a unit, settle the ramified prime-exponent
case, extend the local integer criterion to arbitrary unramified
mixed-characteristic discrete valuation rings, obtain corresponding Dedekind-domain criteria, and classify all
mixed-power difference factor absorbing ideals of $\mathbb Z$.
\end{abstract}
\maketitle

\section{Introduction}

A prime ideal recognizes every product of two elements. A difference-factor
condition recognizes only products arising from a prescribed polynomial
identity. This restriction permits two phenomena that do not occur for
prime ideals: different irreducible components can select different factors,
and nilpotent elements can make the complementary factor vanish without
making the linear factor vanish. Understanding these two mechanisms is
essential for obtaining classifications, rather than only transfer
properties and isolated examples.

Anderson, Badawi, and Coykendall introduced square-difference factor
absorbing ideals in \cite{ABC}. Farshadifar studied the cubic condition in
\cite{Far}, and Ye\c{s}ilot considered prime exponents in \cite{Yes}.
Hizem characterized square-difference factor absorbing ideals using minimal
primes in \cite{Hizem}. The mixed-power formulation of Taha, Salim, and Kim
\cite{TS} starts from the identity
\[
 a^{kn}-b^n=(a^k-b)F_{k,n}(a,b),\qquad
 F_{k,n}(a,b)=\sum_{i=0}^{n-1}a^{k(n-1-i)}b^i.
\]
The companion manuscript establishes the quotient formulation, localization and
product criteria, the complete square case, and several constructions.
It also shows that the multiplier $k$ can matter for nonradical ideals:
$16\mathbb Z$ is $(6,3)$-mpdf but is not $(3,3)$-mpdf.

Tekir, Anebri, Hasi\c{c}, and Arabac{\i}~\cite{TAHA} study a different
quartic condition: from $(a-b)(a+b)(a^2+b^2)\in I$, one requires
one of the three pair products to lie in $I$. We retain their term
\emph{4-pdf-absorbing} for this condition. It is weaker than the
strict $(4,4)$-mpdf condition, which selects $a-b$ or the full cubic
cofactor. Their finite comaximal-prime classification and primary
idealization criterion provide the quartic antecedents for the
comparisons below. Mixed multipliers give further distinctions:
Example~\ref{quartic:mixed-separation} separates 4-pdf absorption
from $(4k,4)$-mpdf for every $k\ge2$.

The purpose of this paper is to identify the structural information that
decides the mixed-power condition. For reduced quotients, that information
is the interaction between roots of unity and residue characteristics at
minimal primes. For nonreduced quotients, it includes power closure,
valuation bounds, and the extension data of square-zero ideals. The
principal contributions are as follows.
\begin{enumerate}
\item Theorem~\ref{thm:spectral} extends the cubic minimal-prime criterion
to arbitrary $n$ and proves independence of $k$, with no Noetherian,
finiteness, or comaximality assumptions.
Theorem~\ref{quartic:radical-criterion} likewise removes finiteness
and comaximality from the quartic pair-product criterion of
\cite[Theorem~2.8]{TAHA}; its obstruction involves three minimal
primes, rather than two.
\item Theorem~\ref{thm:radical-valued} characterizes the implication with
both conclusions in $\sqrt I$ by the mixed-power condition on $\sqrt I$.
It also makes this radical-valued implication independent of $k$ for
arbitrary $I$. Corollary~\ref{thm:radical-descent} yields descent from
the strict condition and from the corresponding primary condition.
Theorem~\ref{thm:series} gives polynomial and formal power series
equivalences for radical ideals. Theorems~\ref{thm:geometric} and
\ref{thm:domain-pullback} give geometric and pullback applications.
\item Theorem~\ref{thm:exact-p-part} shows that, in characteristic $p$ with
$p\mid n$, the condition depends exactly on $p^{v_p(n)}$ and is independent
of $k$ and of the prime-to-$p$ part of $n$. It also yields a local criterion.
\item Theorem~\ref{thm:squarezero-domain} classifies square-zero extensions
of domains without assuming that the extension splits.
Corollary~\ref{cor:arbitrary-idealization-mpdf} includes nonhomogeneous
ideals of an idealization and records the role of their defining linear map.
Theorem~\ref{thm:geometric-squarezero} characterizes stability under
every field extension by geometric integrality of the reduced quotient
and torsion-freeness of the kernel, including nonsplit extensions.
Corollary~\ref{quartic:homogeneous-idealization} treats homogeneous
idealizations over primary bases with regular module action and
removes the restriction on $3$ from \cite[Theorem~5.5]{TAHA}.
\item Theorem~\ref{thm:tame-dvr} gives the exact dependence on $k$ for powers
of the maximal ideal of a discrete valuation ring when $n$ is a unit.
Theorem~\ref{thm:ramified-prime} settles the ramified prime-exponent question
in \cite[Remark~5]{Yes}. Theorem~\ref{thm:unramified-dvr} treats all exponents in unramified
mixed characteristic and identifies the residue-field restriction in the
characteristic-four exception. Theorems~\ref{thm:dedekind-prime} and
\ref{thm:integer-all} give complete classifications in the stated
Dedekind-domain and integer settings.
\end{enumerate}

The square minimal-prime criterion is due to Hizem~\cite{Hizem}.
For cubes, three definitions must be distinguished. Farshadifar's
cdf-absorbing condition requires $a-b\in I$ or
$a^2+ab+b^2\in I$; the cdf-primary condition of
Ye\c{s}ilot, \c{S}im\c{s}ek, and \"Ozel Ay~\cite{YSO} requires
$a-b\in I$ or $a^2+ab+b^2\in\sqrt I$; Draoui's
cdf-semiprimary condition~\cite[Definition~2.1]{Dr26} places both
conclusions in $\sqrt I$. All three agree for radical ideals.
The final published \cite[Theorem~1]{YSO} already gives the cubic
spectral criterion without finiteness or comaximality assumptions,
and \cite[Theorem~2]{YSO} proves cubic primary descent to the
radical in arbitrary rings. We use their two-prime localization and
nil-ideal lifting framework explicitly. The extensions here allow
arbitrary $n$ and mixed multipliers: the order of an individual root
of unity must be separated from $n$, and the vanishing of $n$ at a
second minimal prime is an independent condition.

Theorem~\ref{thm:radical-valued} also treats the weaker implication
with both factors tested in the radical. For nonradical ideals, the
strict condition studied in the remaining classifications retains
information lost by the primary and semiprimary conditions. In
particular, the bounded exponents in
Corollary~\ref{cor:integer-cubic} differ from the unrestricted
prime-power exponents in the cdf-primary classification
\cite[Theorem~6]{YSO}; its squarefree case is already contained in
\cite[Theorem~5]{YSO}.

The prime-exponent case $k=1$ of the tame valuation theorem recovers
\cite[Theorem~3]{Yes}, while Theorem~\ref{thm:ramified-prime} resolves
the question posed in \cite[Remark~5]{Yes}. The latter source is a
preprint. The numerical restrictions for closed prime-power ideals
are already present in Anderson and Badawi~\cite{ABclosed}; we state
their role explicitly rather than treating them as independent
novelty. General homomorphism and product criteria are developed in
the companion manuscript~\cite{TS}. Here they serve as tools for
the spectral, extension and arithmetic classifications. The analogous
homomorphism, quotient, and localization transfers for quartic pair
absorption are \cite[Proposition~3.2, Corollary~3.3, and
Proposition~3.6]{TAHA}. Its product statement
\cite[Proposition~3.1]{TAHA} concerns ideals with one whole-ring
coordinate; the strict compatibility criterion used here comes from
\cite{TS}.

All conditions in this paper quantify over every pair $a,b$, including
zero. This agrees with the published cubic definition in
\cite[Definition~2.1]{Far} and
with \cite{Yes,TS}. The original square definition in \cite{ABC} excludes
zero entries. The two square conventions agree for nonzero ideals, but
need not agree for the zero ideal; the quotient-based statements below
always use the all-elements convention.

Section~\ref{sec:prelim} collects the transfer and localization tools.
Section~\ref{sec:spectral} gives the strict and quartic pair-product
spectral criteria, radical descent, and extension applications.
Section~\ref{sec:charp} treats the characteristic reduction, while
Section~\ref{sec:squarezero} treats square-zero extensions and the
homogeneous-primary quartic comparison.
Sections~\ref{sec:valuation} and~\ref{sec:integers} give the valuation,
Dedekind-domain, and integer classifications.

\section{Definition, reformulations and first examples}\label{sec:prelim}

All rings are commutative with identity. Except for an auxiliary target
of a pullback, rings are nonzero. Homomorphisms preserve the identity.
An ideal under consideration is proper unless stated otherwise. We fix
$n\ge2$ and $k\ge1$, and write
$F_n(X,Y)=\sum_{i=0}^{n-1}X^{n-1-i}Y^i$, so
$F_{k,n}(a,b)=F_n(a^k,b)$. We put $F_1=1$. For a field $K$, write
$\mu_n(K)=\{z\in K:z^n=1\}$. The notation $n\in I$ means
$n\cdot1_R\in I$. General facts about radicals, localization, and the
Chinese remainder theorem are used as in \cite{AM}.

\begin{definition} \cite{TS}
A proper ideal $I$ of $R$ is \emph{$(kn,n)$-mixed-power difference factor
absorbing}, abbreviated \emph{$(kn,n)$-mpdf}, if
$a^{kn}-b^n\in I$ implies $a^k-b\in I$ or $F_{k,n}(a,b)\in I$ for all
$a,b\in R$. A nonzero ring $S$ satisfies $\MP{k}{n}$ if its zero ideal
is $(kn,n)$-mpdf.
\end{definition}

\begin{definition}\label{quartic:definition}
A proper ideal $I$ of $R$ is \emph{4-pdf-absorbing}
\cite[Definition~2.1]{TAHA} if, for all $a,b\in R$,
\[
 a^4-b^4\in I\quad\Longrightarrow\quad
 \begin{aligned}[t]
 &a^2-b^2\in I,\quad\text{or}\\
 &(a^2+b^2)(a-b)\in I,\quad\text{or}\\
 &(a^2+b^2)(a+b)\in I.
 \end{aligned}
\]
We also call this the \emph{quartic pair-product condition}.
\end{definition}

Every $(4,4)$-mpdf ideal is 4-pdf-absorbing: membership of $a-b$
implies membership of $a^2-b^2$, and the other strict alternative
is already a pair product. The converse fails, as illustrated in
Example~\ref{quartic:radical-examples}. Both conditions are tested
in the quotient $R/I$.

\begin{example}
Let $R=\mathbb F_2[x]$, and let $I=(x^2+x)=(x(x+1))$.
We show that $I$ is $(2k,2)$-mpdf for every $k\ge1$.
It is not prime, since $x(x+1)\in I$ while $x,x+1\notin I$.

Let $a,b\in R$ satisfy $a^{2k}-b^2\in I$, and put $c=a^k$.
Since $\operatorname{char}R=2$, we have $-1=1$ and $2=0$, so
\[
(c-b)^2=c^2-2cb+b^2=c^2+b^2=c^2-b^2=a^{2k}-b^2\in I.
\]
Hence $x(x+1)$ divides $(c-b)^2$.
The elements $x$ and $x+1$ are prime in $R$, so $x\mid c-b$ and $x+1\mid c-b$.
They are coprime, so $x(x+1)\mid c-b$, that is, $a^k-b\in I$.
\end{example}

We recall the following mixed-power transfer statements from \cite{TS}.
For $n=3$ and $k=1$, the power-closure implication is
\cite[Proposition~2.4]{Far}, localization is
\cite[Theorem~2.15]{Far}, and the homomorphism and quotient
statements are \cite[Theorem~2.16 and Corollary~2.17]{Far}.

\begin{lemma}\label{lem:transfers}
Let $I$ be a proper ideal of $R$.
\begin{enumerate}
\item $I$ is $(kn,n)$-mpdf if and only if $R/I$ satisfies $\MP{k}{n}$.
\item The property $\MP{k}{n}$ passes to unital subrings and to every
nonzero localization. It therefore passes both ways between a ring and
its total quotient ring.
\item If $I$ is $(kn,n)$-mpdf, then $c^n\in I$ implies $c^{n-1}\in I$.
\end{enumerate}
\end{lemma}

For the quartic pair-product condition, the corresponding localization
statement is \cite[Proposition~3.6]{TAHA}. We record the
denominator-clearing argument used below.

\begin{lemma}\label{quartic:localization}
If the zero ideal of a nonzero ring $S$ is 4-pdf-absorbing, then the
zero ideal of every nonzero localization of $S$ is 4-pdf-absorbing.
\end{lemma}
\begin{proof}
If $(a/s)^4=(b/t)^4$ in a localization, choose $q$ in its
multiplicative set such that $q(a^4t^4-b^4s^4)=0$.
Put $A=qat$, $B=qbs$, and $d=qst$. Then $A^4=B^4$, so one
of the three pair products of $A,B$ is zero. Since $a/s=A/d$
and $b/t=B/d$, division by $d^2$ or $d^3$ gives the corresponding
conclusion in the localization.
\end{proof}

The two-minimal-prime localization used in
\cite[proof of Theorem~1]{YSO} has the following finite-family form.
The cases of two and three primes will serve the strict and quartic
pair-product conditions, respectively.

\begin{lemma}\label{lem:two-points}\label{quartic:finite-points}
Let $P_1,\ldots,P_r$ be distinct minimal prime ideals of a nonzero
ring $R$, where $r\ge1$, and put
$T=R\setminus\bigcup_{j=1}^rP_j$. The ring $C=T^{-1}R$ has
exactly $r$ prime ideals, $T^{-1}P_1,\ldots,T^{-1}P_r$, all maximal.
In particular,
\[
 C/\sqrt{0_C}\cong\prod_{j=1}^rK_{P_j},
 \qquad K_{P_j}=\Frac(R/P_j).
\]
If $R$ is reduced, then $C\cong\prod_{j=1}^rK_{P_j}$.
\end{lemma}
\begin{proof}
The set $T$ is multiplicatively closed. A prime $H$ disjoint from
$T$ is contained in $\bigcup_{j=1}^rP_j$. Finite prime avoidance
gives $H\subseteq P_j$ for some $j$, and minimality forces
$H=P_j$. Conversely, every $P_j$ is disjoint from $T$.
Thus exactly these $r$ primes survive; all are maximal and hence
pairwise comaximal. Their intersection is the nilradical, so the
Chinese remainder theorem gives the displayed product. Each
$C/P_jC$ is a field localization of $R/P_j$, hence is its fraction
field. Finally, localization preserves reducedness.
\end{proof}

\section{Radical ideals and their minimal primes}\label{sec:spectral}

In a field, equality of the two powers is controlled by roots of unity.
The next theorem identifies precisely when the choices at different
minimal primes are compatible. For $n=3$ and $k=1$, this is the
unrestricted cubic criterion of \cite[Theorem~1]{YSO}.

\begin{theorem}\label{thm:spectral}
Let $I$ be a proper radical ideal of a commutative ring $R$, let
$n\ge2$, and put $K_P=\operatorname{Frac}(R/P)$ for
$P\in\operatorname{Min}_R(I)$.  The following are equivalent:
\begin{enumerate}
\item $I$ is $(kn,n)$-mpdf for some $k\ge1$;
\item $I$ is $(kn,n)$-mpdf for every $k\ge1$;
\item for all distinct $P,Q\in\operatorname{Min}_R(I)$,
\[
 \mu_n(K_P)\ne\{1\}\quad\text{then}\quad n\in Q;
\]
\item for every distinct $P,Q\in\operatorname{Min}_R(I)$, the zero
ideal of $K_P\times K_Q$ is $(kn,n)$-mpdf (for one, equivalently all,
positive integers $k$).
\end{enumerate}
\end{theorem}
\begin{proof}
Pass to the reduced quotient $S=R/I$.  For a field $K$ and
$a^{kn}=b^n$, either $a=b=0$, or $a\ne0$ and
$\zeta=b/a^k\in\mu_n(K)$.  The first factor is nonzero exactly when
$\zeta\ne1$.  When the first factor is zero, the cofactor is
$na^{k(n-1)}$; when the first factor is nonzero, the cofactor vanishes.

Suppose $S$ satisfies $\MP{k}{n}$.  Fix distinct minimal primes
$P,Q$.  Lemma~\ref{lem:two-points} and localization transfer imply
that $K_P\times K_Q$ satisfies $\MP{k}{n}$.  If
$\zeta\in\mu_n(K_P)\setminus\{1\}$ and $n\ne0$ in $K_Q$, use
$a=(1,1)$ and $b=(\zeta,1)$.  Then $a^{kn}=b^n$, whereas
$a^k-b=(1-\zeta,0)\ne0$ and the cofactor has second coordinate
$n\ne0$.  This is impossible.

Conversely, suppose the condition on distinct minimal primes holds.
Let $a^{kn}=b^n$ in $S$, and suppose the first factor $a^k-b$ is
nonzero.  Since the intersection of all minimal primes is zero, some
minimal prime $P$ does not contain $a^k-b$.  In $K_P$ the ratio
$b/a^k$ is a nontrivial $n$th root of unity.  Thus $n\in Q$ for every
other minimal prime $Q$.  Modulo $P$ the cofactor vanishes by the
domain property.  Modulo any other $Q$, either the first factor is
nonzero, in which case the cofactor again vanishes, or $a^k=b$, in
which case the cofactor equals $na^{k(n-1)}=0$.  Therefore the cofactor
belongs to every minimal prime and is zero.  This works for every
$k\ge1$.

For (iii) implies (iv), fix distinct $P,Q$ and put $W=K_P\times K_Q$.
The minimal primes of $W$ are $0\times K_Q$ and $K_P\times0$, with residue fields
$K_P$ and $K_Q$. Condition (iii) for $W$ consists of the two instances $(P,Q)$ and
$(Q,P)$ of condition (iii) for $S$. The argument just given, applied to $W$,
shows that $W$ satisfies $\MP{k}{n}$ for every $k$.
Conversely, the computation with $a=(1,1)$ and $b=(\zeta,1)$ in $K_P\times K_Q$,
used above to show that $\MP{k}{n}$ for $S$ implies (iii), shows that (iv) implies (iii).
\end{proof}

\begin{remark}
For a reduced ring $S$, set
$\mathcal A_n=\{P\in\Min S:\mu_n(K_P)\ne\{1\}\}$ and
$\mathcal B_n=\{P\in\Min S:n\notin P\}$. The criterion says that
$\mathcal A_n\times\mathcal B_n$ is contained in the diagonal.
Equivalently, either $\mathcal A_n=\varnothing$, or
$\mathcal B_n=\varnothing$, or
$\mathcal A_n=\mathcal B_n=\{P_0\}$ for one minimal prime $P_0$.
The first alternative forces the $n$th-power map on $S$ to be
injective.  The second means $n=0$ in the reduced ring $S$ and ensures
that the cofactor always vanishes whenever the mixed powers agree.
These alternatives are not necessarily mutually exclusive.
\end{remark}

\begin{corollary}\label{cor:square-recovery}
Let $I$ be a proper ideal of $R$ and let $k\ge1$. Then $I$ is
$(2k,2)$-mpdf if and only if it is radical and at most one prime minimal
over $I$ fails to contain $2$.
\end{corollary}
\begin{proof}
Lemma~\ref{lem:transfers}(iii) gives $c^2\in I\Rightarrow c\in I$.
If $c^m\in I$ with $m\ge2$, then $(c^{\lceil m/2\rceil})^2\in I$, so $c^{\lceil m/2\rceil}\in I$,
and repeating this lowers the exponent until $c\in I$. Hence a $(2k,2)$-mpdf ideal is radical.
For a field $K$, the group $\mu_2(K)$ is nontrivial exactly when
$2\ne0$ in $K$. Apply Theorem~\ref{thm:spectral}.
\end{proof}

This is the square-case theorem of \cite{TS}, extending the
minimal-prime description of \cite{Hizem}, with the all-elements convention for the zero ideal.
For general $n$, the two sets in the preceding remark can be different.
For example, with $n=6$ and $K=\mathbb F_4$, one has $6=0$ in $K$
but $\mu_6(K)\ne\{1\}$: a generator of $K^\times$ has order $3$.
Thus nontrivial $n$th roots of unity can remain in a characteristic
dividing $n$.

\subsection{Quartic pair absorption}
\label{quartic:spectral-section}

The quartic pair-product condition has a different compatibility
criterion. For an intersection of finitely many pairwise comaximal
primes, the following classification is
\cite[Theorem~2.8]{TAHA}. Localization at three selected minimal
primes removes both restrictions.

\begin{theorem}\label{quartic:radical-criterion}
Let $I$ be a proper radical ideal of $R$. For $P\in\Min_R(I)$ put
$K_P=\Frac(R/P)$, and define
\[
 \Sigma=\{P\in\Min_R(I):2\notin P\},\qquad
 \Sigma_i=\{P\in\Sigma:-1\text{ is a square in }K_P\}.
\]
Then
\[
 I\text{ is 4-pdf-absorbing}
 \quad\Longleftrightarrow\quad
 \Sigma_i=\varnothing\ \text{or}\ |\Sigma|\le2.
\]
\end{theorem}
\begin{proof}
Pass to the reduced quotient $S=R/I$, identifying its minimal
primes and their fraction fields with those in the statement.
Suppose $a^4=b^4$, and put $u=a-b$, $v=a+b$, $w=a^2+b^2$.
At a minimal prime containing $2$, the residue domain has
characteristic two. There $(a-b)^4=0$ gives $a=b$, and all
three factors $u,v,w$ vanish.

Suppose that none of $uv,uw,vw$ is zero in $S$. Since the
intersection of its minimal primes is zero, choose minimal
primes $P,Q,H$ at which $uv,uw,vw$, respectively, are nonzero.
At $P$, the identity $uvw=0$ forces $w=0$ while $u,v\ne0$.
Thus $a,b\ne0$ and $(a/b)^2=-1$ in $K_P$, so $P\in\Sigma_i$.
At $Q$ the same identity forces $v=0$, while $u,w\ne0$;
at $H$ it forces $u=0$, while $v,w\ne0$.
All three primes belong to $\Sigma$, and they are distinct:
at each of them the other two pair products vanish.
Hence failure of the 4-pdf condition implies
$\Sigma_i\ne\varnothing$ and $|\Sigma|\ge3$. This proves
sufficiency of the displayed criterion.

Conversely, suppose $P\in\Sigma_i$ and choose distinct
$Q,H\in\Sigma\setminus\{P\}$. By
Lemma~\ref{quartic:finite-points}, localization away from
$P\cup Q\cup H$ gives $K_P\times K_Q\times K_H$.
Choose $i\in K_P$ with $i^2=-1$, and take
\[
 a=(1,1,1),\qquad b=(i,-1,1).
\]
Then $a^4=b^4$, whereas
\[
 \begin{aligned}
 a^2-b^2&=(2,0,0),\\
 (a^2+b^2)(a-b)&=(0,4,0),\\
 (a^2+b^2)(a+b)&=(0,0,4).
 \end{aligned}
\]
All three are nonzero since these fields have characteristic
different from two. Lemma~\ref{quartic:localization} therefore
precludes the 4-pdf condition in $S$, proving necessity.
\end{proof}

\begin{remark}\label{quartic:scope-of-extension}
If $\Sigma_i=\varnothing$, the first conclusion always suffices:
$a^4-b^4\in I$ implies $a^2-b^2\in I$. Indeed, in a residue
domain of characteristic different from two, $w=0$ with $b\ne0$
would produce a square root of $-1$; otherwise $uv=0$.
In characteristic two all three factors vanish as in the proof.
Thus the finiteness and comaximality assumptions in
\cite[Corollary~2.11]{TAHA} may be removed. Over an algebraically
closed field of characteristic different from two, every
minimal-prime fraction field contains a square root of $-1$.
Consequently, in any algebra over that field, a proper radical
ideal is 4-pdf-absorbing exactly when it has at most two minimal
primes. This extends \cite[Corollary~2.12]{TAHA}, whose statement
concerns a polynomial ring and assumes comaximality.
\end{remark}

\begin{corollary}\label{quartic:strict-radical-comparison}
With the notation of Theorem~\ref{quartic:radical-criterion}, the
radical ideal $I$ is $(4k,4)$-mpdf for some, equivalently every,
$k\ge1$ if and only if $|\Sigma|\le1$.
\end{corollary}
\begin{proof}
For a field $K$, $\mu_4(K)\ne\{1\}$ exactly when its
characteristic is different from two: then $-1\ne1$ is a fourth
root, whereas in characteristic two $X^4-1=(X-1)^4$.
Also $4\in P$ exactly when $2\in P$. The criterion of
Theorem~\ref{thm:spectral} therefore forbids precisely two
distinct members of $\Sigma$.
\end{proof}

\begin{example}\label{quartic:radical-examples}
The following examples separate the two quartic conditions.
\begin{enumerate}
\item The zero ideal of $\Q^r$ is 4-pdf-absorbing for every
positive integer $r$, because $-1$ is not a square in $\Q$.
It is $(4k,4)$-mpdf only when $r=1$. For $r\ge2$, take
$a=(1,\ldots,1)$ and $b=(1,-1,1,\ldots,1)$: the linear
factor is nonzero in the second coordinate and the cofactor
is $4$ in the first, for every $k\ge1$.
\item The zero ideal of $\C^2$ is 4-pdf-absorbing but is not
$(4k,4)$-mpdf, by the same two-coordinate pair.
The zero ideal of $\C^3$ is not 4-pdf-absorbing, by the witness
in Theorem~\ref{quartic:radical-criterion}. These product examples
already fall within \cite[Theorem~2.8]{TAHA}.
\item In $S=\C[x,y]/(xy)$, the two minimal primes $(x)$ and
$(y)$ are not comaximal. Its zero ideal is nevertheless
4-pdf-absorbing by Theorem~\ref{quartic:radical-criterion}.
For the residue classes $a=x+y$, $b=x-y$ one has $a^4=b^4$,
but $a-b=2y\ne0$ and $(a+b)(a^2+b^2)=4x^3\ne0$.
Thus the zero ideal is not $(4,4)$-mpdf. This example lies
outside the comaximal hypothesis in \cite[Theorem~2.8]{TAHA}.
\end{enumerate}
\end{example}

\subsection{Descent to radicals}

The cubic primary descent theorem \cite[Theorem~2]{YSO} uses
two-prime localization and lifting across the nilradical. We extend
this mechanism to arbitrary exponents and to the weaker implication
with both conclusions in the radical. The standard lifting facts
for idempotents and simple roots across nil ideals appear in
\cite[Tags 00J9 and 0ALI]{Stacks}. We give elementary
proofs of the precise forms needed, so no nilpotence assumption on the
entire nilradical is required.

\begin{lemma}\label{lem:nil-lifts1} \cite[Theorem~21.28]{Lam2001}
Let $N$ be a nil ideal of a commutative ring $C$. Every idempotent of $C/N$ lifts
to an idempotent of $C$.
\end{lemma}

\begin{lemma}\label{lem:nil-lifts2}
Let $N$ be a nil ideal of a commutative ring $C$, and let $d\ge1$ be a unit in $C$.
Every $\bar z\in C/N$ with $\bar z^d=1$ lifts to $z\in C$ with $z^d=1$.
\end{lemma}

\begin{proof}
We use the following elementary fact. If $g\in C$ is nilpotent, say
$g^m=0$, then $1+g$ is a unit, since
\[
 (1+g)\bigl(1-g+g^2-\cdots+(-g)^{m-1}\bigr)=1-(-g)^m=1.
\]

Let $\bar z\in C/N$ satisfy $\bar z^d=1$, and choose $z_0\in C$ with
$z_0+N=\bar z$. Then $h=z_0^d-1$ lies in $N$, so $h^m=0$ for some
$m\ge1$. The element $z_0^d=1+h$ is a unit, and since
$z_0\cdot z_0^{d-1}=z_0^d$, the element $z_0$ is a unit as well. We construct
$z_0,z_1,z_2,\ldots$ in $C$ such that, for every $j\ge0$,
\[
 z_j\equiv z_0 \pmod N,\qquad z_j^d-1\in h^{2^j}C.
\]
For $j=0$ there is nothing to prove. Suppose $z=z_j$ has been
constructed, and put $g=z^d-1\in h^{2^j}C\subseteq N$. Since
$z\equiv z_0\pmod N$, we have
\[
 z=z_0\bigl(1+z_0^{-1}(z-z_0)\bigr),\qquad z_0^{-1}(z-z_0)\in N,
\]
so $z$ is a unit by the fact recorded above. Since $d$ is a unit as well,
\[
 w=d\,z^{d-1}
\]
is a unit of $C$. Put $t=-g\,w^{-1}\in N$ and $z'=z+t$. By the binomial
theorem,
\[
 z'^d=z^d+d\,z^{d-1}t+t^2q,\qquad
 q=\sum_{i=2}^{d}\binom di z^{d-i}t^{i-2}\in C,
\]
and therefore
\[
 z'^d-1=g+wt+t^2q=g-g+g^2w^{-2}q=g^2w^{-2}q\in g^2C.
\]
Since $g^2\in h^{2^{j+1}}C$, we get $z'^d-1\in h^{2^{j+1}}C$. Moreover
$z'\equiv z\equiv z_0\pmod N$, because $t\in N$. We take $z_{j+1}=z'$.

Finally, choose $j$ with $2^j\ge m$. Then $h^{2^j}=0$, so $z_j^d=1$, and
$z_j+N=\bar z$. Thus $z=z_j$ is the required lift.
\end{proof}

The next theorem determines exactly when testing the two factors in the
radical suffices. Its hypothesis is weaker than the strict mixed-power
condition and than the corresponding primary condition.

\begin{theorem}\label{thm:radical-valued}
Let $R$ be a commutative ring, let $I\subsetneq R$, and fix integers
$n\ge2$ and $k\ge1$. The following conditions are equivalent.
\begin{enumerate}
\item For every $a,b\in R$,
\[
 a^{kn}-b^n\in I\quad\Longrightarrow\quad
 a^k-b\in\sqrt I\ \text{or}\ F_{k,n}(a,b)\in\sqrt I.
\]
\item The radical $\sqrt I$ is $(kn,n)$-mpdf.
\item The radical $\sqrt I$ is $(\ell n,n)$-mpdf for every $\ell\ge1$.
\end{enumerate}
Consequently, for $n=3$ and $k=1$, an ideal is cdf-absorbing semiprimary
in the sense of \cite[Definition~2.1]{Dr26} if and only if its radical
is cdf-absorbing. For every fixed $n$, condition (i) is independent of $k$.
\end{theorem}
\begin{proof}
Condition (iii) implies (ii), and (ii) implies (i), since $I\subseteq\sqrt I$.
Assume (i), and put $S=R/I$ and $N=\sqrt{0_S}$. The hypothesis says
that whenever $a^{kn}=b^n$ in $S$, at least one of $a^k-b$ and
$F_{k,n}(a,b)$ is nilpotent.

We first verify that this weaker implication passes to every nonzero
localization $T^{-1}S$. For $n=3$ and $k=1$, this localization property
is \cite[Theorem~5.3]{Dr26}; we verify the weighted version needed here.
Suppose $A=a/s$, $B=b/t$ satisfy
$A^{kn}=B^n$. Put $w=st$, $a'=at$, and $b'=s^kt^{k-1}b$, so that
$a'/1=wA$ and $b'/1=w^kB$. There is $c\in T$ such that
$c(a'^{kn}-b'^n)=0$. Consequently,
$(ca')^{kn}-(c^kb')^n=0$ in $S$. The hypothesis makes at least one of
the following elements nilpotent:
\[
 c^k(a'^k-b'),\qquad
 c^{k(n-1)}F_{k,n}(a',b')
\]
In the localization these are respectively
$(cw)^k(A^k-B)$ and $(cw)^{k(n-1)}F_{k,n}(A,B)$.
Since $cw$ is a unit, one of the two desired factors is nilpotent.

We now prove the minimal-prime condition of Theorem~\ref{thm:spectral}
for the reduced ring $S/N$, writing $K_P=\operatorname{Frac}(S/P)$
for each $P\in\Min S$. Suppose instead that distinct minimal
primes $P,Q$ of $S$ satisfy
$\mu_n(K_P)\ne\{1\}$ and $n\notin Q$. Choose
$\zeta\in\mu_n(K_P)\setminus\{1\}$, and let $d>1$ be its order.
Then $d\mid n$, and $d$ is nonzero in $K_P$: a field of
characteristic $p$ has no nontrivial root of unity of order divisible
by $p$. This uses $d$, not $n$; the characteristic of $K_P$ may
divide $n$.

Set $T=S\setminus(P\cup Q)$ and $C=T^{-1}S$. By
Lemma~\ref{lem:two-points},
$C/\sqrt{0_C}\cong K_P\times K_Q$. Lift the idempotent $(1,0)$
using Lemma~\ref{lem:nil-lifts1}, obtaining a decomposition
$C\cong A\times B$ with
$A/\sqrt{0_A}\cong K_P$ and $B/\sqrt{0_B}\cong K_Q$.
The integer $d$ is a unit in $A$, since its image is a unit modulo
a nil ideal. Lemma~\ref{lem:nil-lifts2} therefore lifts $\zeta$ to
$z\in A$ satisfying $z^d=1$ and hence $z^n=1$.

Take $a=(1_A,1_B)$ and $b=(z,1_B)$. Then $a^{kn}=b^n$, while
the reductions of the two factors modulo $\sqrt{0_C}$ are
\[
 (1-\zeta,0)\ne0,
 \qquad
 \left(\sum_{i=0}^{n-1}\zeta^i,\ n\cdot1_{K_Q}\right)\ne0.
\]
Both factors are therefore nonnilpotent. This contradicts the
localized implication already proved. The minimal-prime condition
holds, and Theorem~\ref{thm:spectral} makes $S/N$ satisfy every
mixed multiplier with exponent $n$. Since $S/N\cong R/\sqrt I$,
Lemma~\ref{lem:transfers}(i) gives (iii).
\end{proof}

\begin{corollary}\label{thm:radical-descent}
Let $R$ be a commutative ring, $I$ a proper ideal, $n\ge2$, and
$k\ge1$. If $I$ is $(kn,n)$-mpdf, then $\sqrt I$ is
$(\ell n,n)$-mpdf for every $\ell\ge1$. More generally, the same
conclusion follows from
\[
 a^{kn}-b^n\in I\quad\Longrightarrow\quad
 a^k-b\in I\ \text{or}\ F_{k,n}(a,b)\in\sqrt I
 \qquad(a,b\in R).
\]
\end{corollary}
\begin{proof}
Either hypothesis implies condition (i) of
Theorem~\ref{thm:radical-valued}.
\end{proof}

For $n=3$ and $k=1$, the primary assertion of the corollary is
\cite[Theorem~2]{YSO}. The equivalence in
Theorem~\ref{thm:radical-valued} also corrects the negative assertion
about semiprimary radical descent in \cite[Remark~2.5]{Dr26}:
the examples cited there are explicitly non-semiprimary and do not
provide counterexamples to the equivalence.

\subsection{Polynomial and formal power series extensions}

The absence of new roots of unity in a Laurent series field makes the
spectral classification stable under adjoining variables.

\begin{lemma}\label{lem:Laurent-torsion}
For every field $K$ and every $n\ge2$,
$\mu_n(K((X)))=\mu_n(K)$.
\end{lemma}
\begin{proof}
A torsion Laurent series has valuation zero.  Its constant coefficient
$c$ satisfies $c^n=1$.  Divide by $c$, obtaining $1+w$ with
$w\in XK[[X]]$.  In characteristic zero set $r=n$; in characteristic
$p$ write $n=p^er$ with $p\nmid r$.  Since Frobenius is injective in
the field $K((X))$, the equation $(1+w)^n=1$ implies
$(1+w)^r=1$.  If $w$ is nonzero of order $s\ge1$, the coefficient
of $X^s$ in $(1+w)^r-1$ is $r$ times the nonzero leading coefficient
of $w$, a contradiction.  Hence $w=0$.
\end{proof}

The cubic polynomial constructions in \cite[Theorem~5.1]{Dr26} and
\cite[Theorem~7(2)]{YSO} concern $(I,X)$ and preserve the quotient
$R/I$. For quartic pair absorption, \cite[Proposition~3.9]{TAHA}
proves descent from $I[X]$ to $I$, while \cite[Example~3.10]{TAHA}
shows failure of general ascent from $\Z/8\Z$ to its polynomial
ring. The next theorem establishes equivalence for the strict
condition on radical coefficient ideals $I[X]$ and $I[[X]]$.
Recall that $I[[X_1,\ldots,X_d]]$ means all series whose coefficients lie in $I$.  It equals $IR[[X_1,\ldots,X_d]]$ when $I$ is finitely
generated, but this equality does not hold in general.

\begin{theorem}\label{thm:series}
Let $I$ be a proper radical ideal of a commutative ring $R$, let
$n\ge2$, $k\ge1$, and $d\ge1$.  The following are equivalent:
\begin{enumerate}
\item $I$ is $(kn,n)$-mpdf in $R$;
\item $I[X_1,\ldots,X_d]$ is $(kn,n)$-mpdf in
$R[X_1,\ldots,X_d]$;
\item the coefficient ideal $I[[X_1,\ldots,X_d]]$ is
$(kn,n)$-mpdf in $R[[X_1,\ldots,X_d]]$.
\end{enumerate}
\end{theorem}
\begin{proof}
Pass to $S=R/I$.  The constant embeddings give both reverse
implications.  Suppose $S$ is reduced and satisfies
$\MP{k}{n}$.  The coefficientwise maps give an embedding
\[
 S[[X]]\ \longrightarrow\
 \prod_{P\in\operatorname{Min}S}K_P[[X]]\ \longrightarrow\
 \prod_{P\in\operatorname{Min}S}K_P((X)).
\]
By Lemma~\ref{lem:Laurent-torsion}, the fields in the last product
have exactly the same $n$th roots of unity and characteristics as
the corresponding $K_P$.  Write $L_P=K_P((X))$, and let $a,b\in S[[X]]$ satisfy $a^{kn}=b^n$ and $a^k\ne b$.
Choose $P\in\operatorname{Min}S$ whose coordinate satisfies $a_P^k\ne b_P$ in $L_P$.
Then $a_P\ne0$, since otherwise $b_P^n=0$ forces $b_P=0$, and
$\zeta=b_P/a_P^k\in\mu_n(L_P)\setminus\{1\}$.
By Lemma~\ref{lem:Laurent-torsion}, $\mu_n(K_P)\ne\{1\}$, so
Theorem~\ref{thm:spectral}(iii) gives $n\in Q$ for every minimal prime $Q\ne P$.
As $\operatorname{char}L_Q=\operatorname{char}K_Q$, we get $n=0$ in $L_Q$ for $Q\ne P$.
Now $F_{k,n}(a_P,b_P)=0$, because $L_P$ is a domain and $\zeta\ne1$.
For $Q\ne P$, either $a_Q^k\ne b_Q$, and then $F_{k,n}(a_Q,b_Q)=0$ in the field $L_Q$,
or $a_Q^k=b_Q$, and then $F_{k,n}(a_Q,b_Q)=n\,a_Q^{k(n-1)}=0$.
Thus every coordinate of $F_{k,n}(a,b)$ vanishes, and $F_{k,n}(a,b)=0$ because
$S[[X]]\to\prod_PL_P$ is injective.
So $S[[X]]$ satisfies $\MP{k}{n}$, and so does its subring $S[X]$.
Both are reduced, so iteration proves the assertion for finitely
many variables.  The quotient identities for coefficient ideals
finish the proof.
\end{proof}

\subsection{Field extension and geometric integrality}

A nontrivial root of unity already present in the base field constrains
every irreducible component simultaneously.

\begin{corollary}\label{cor:base-field}
Let $F$ be a field, $n\ge2$, $k\ge1$, and suppose $\operatorname{char}F$
does not divide $n$ and $\mu_n(F)\ne\{1\}$.  In every $F$-algebra,
the radical $(kn,n)$-mpdf ideals are exactly the prime ideals.
Over an algebraically closed field $F$ with
$\operatorname{char}F\nmid n$ the same conclusion holds.
If $\operatorname{char}F\mid n$, then over any field $F$ every proper radical ideal of
every $F$-algebra is $(kn,n)$-mpdf for every $k\ge1$.
\end{corollary}
\begin{proof}
Every minimal fraction field contains $F$, so under the first
hypotheses it contains a nontrivial $n$th root of unity, and $n$ lies
in no minimal prime.
The spectral criterion allows at most one minimal prime; a radical
ideal with a unique minimal prime is prime.  If $F$ is algebraically
closed and its characteristic does not divide $n$, the polynomial
$X^n-1$ has $n\ge2$ distinct roots.  If the characteristic divides
$n$, the integer $n$ is zero in every quotient and the spectral
criterion is automatic.
\end{proof}

An $F$-algebra $S$ is \emph{geometrically reduced}, respectively
\emph{geometrically integral}, if $S\otimes_F L$ is reduced,
respectively a domain, for every field extension $L/F$.
We shall use the following standard fact about tensor products over an algebraically closed field.

\begin{lemma}\label{lem:tensor-domain}
Let $K$ be an algebraically closed field, and let $D$ and $B$ be $K$-algebras that are integral domains. Then $D\otimes_K B$ is an integral domain.
\end{lemma}
\begin{proof}
The ring $D\otimes_K B$ is nonzero, because $D$ and $B$ are nonzero $K$-vector spaces.
Let $f,g\in D\otimes_K B$ satisfy $fg=0$, and suppose, to the contrary, that $f\neq0$ and $g\neq0$.
Write
\[ f=\sum_{i=1}^{r}a_i\otimes b_i,\qquad g=\sum_{j=1}^{s}c_j\otimes d_j, \]
with $a_i,c_j\in D$, where $b_1,\dots,b_r$ are linearly independent over $K$ and so are $d_1,\dots,d_s$; this is possible after choosing a $K$-basis of the span of the elements of $B$ that occur.
Since $f\neq0$ and $g\neq0$, there are indices $p,q$ with $a_p\neq0$ and $c_q\neq0$.
Let $A$ be the $K$-subalgebra of $D$ generated by all $a_i$ and $c_j$.
Then $A$ is a finitely generated $K$-algebra and a domain, and $A\otimes_K B\to D\otimes_K B$ is injective because $B$ is flat over the field $K$.
Hence we may regard $f$ and $g$ as elements of $A\otimes_K B$ with $fg=0$.

Since $A$ is a domain, $a_pc_q\neq0$, so $a_pc_q$ is not nilpotent and the localization $A'=A[1/(a_pc_q)]$ is a nonzero finitely generated $K$-algebra.
Let $M$ be a maximal ideal of $A'$.
The field $A'/M$ is a finitely generated $K$-algebra, hence a finite extension of $K$ by Zariski's lemma \cite[Proposition~7.9]{AM}, and therefore equal to $K$, because $K$ is algebraically closed.
The composite $A\to A'\to A'/M=K$ is thus a $K$-algebra homomorphism $x:A\to K$.
It sends $a_pc_q$ to the image of a unit of $A'$, so $x(a_p)x(c_q)\neq0$.

Now consider the ring homomorphism $\varphi=x\otimes\mathrm{id}_B:A\otimes_K B\to K\otimes_K B=B$.
We have
\[ \varphi(f)=\sum_{i=1}^{r}x(a_i)\,b_i,\qquad \varphi(g)=\sum_{j=1}^{s}x(c_j)\,d_j. \]
Since the $b_i$ are linearly independent over $K$ and $x(a_p)\neq0$, we get $\varphi(f)\neq0$; likewise $\varphi(g)\neq0$.
But $\varphi(f)\varphi(g)=\varphi(fg)=0$ in the domain $B$, a contradiction.
Hence $f=0$ or $g=0$.
\end{proof}

The next statement identifies the mpdf condition that remains valid
after arbitrary extension of the base field.

\begin{theorem}\label{thm:geometric}
Let $F$ be a field, let $S$ be a nonzero geometrically reduced
$F$-algebra, and let $n\ge2$ with
$\operatorname{char}F\nmid n$.  For any $k\ge1$, the following
are equivalent:
\begin{enumerate}
\item $S\otimes_F L$ satisfies $\MP{k}{n}$ for every field
extension $L/F$;
\item $S\otimes_F\overline F$ satisfies $\MP{k}{n}$;
\item $S$ is geometrically integral over $F$.
\end{enumerate}
\end{theorem}
\begin{proof}
Assume (i).
Taking $L=\overline F$ gives (ii).

Assume (ii), and put $T=S\otimes_F\overline F$.
Since $S$ and $\overline F$ are nonzero $F$-vector spaces, $T\neq0$, and $T$ is reduced because $S$ is geometrically reduced.
The ring $T$ is an algebra over the algebraically closed field $\overline F$, whose characteristic does not divide $n$.
By (ii), the zero ideal of $T$ is $(kn,n)$-mpdf, and it is a proper radical ideal.
By Corollary~\ref{cor:base-field}, applied with $\overline F$ in place of $F$, the zero ideal of $T$ is prime, that is, $T$ is a domain.
Now let $L/F$ be any field extension, and let $\Omega$ be an algebraic closure of $L$.
Since $\overline F$ is algebraic over $F$, there is an $F$-embedding $\overline F\to\Omega$, and through it $\Omega$ is an extension field of $\overline F$.
Then
\[ T\otimes_{\overline F}\Omega=(S\otimes_F\overline F)\otimes_{\overline F}\Omega\cong S\otimes_F\Omega, \]
and this ring is a domain by Lemma~\ref{lem:tensor-domain}.
Since $S$ is a free $F$-module, the inclusion $L\subseteq\Omega$ induces an injective ring homomorphism $S\otimes_F L\to S\otimes_F\Omega$.
Hence $S\otimes_F L$ is a nonzero subring of a domain, so it is a domain.
As $L$ was arbitrary, $S$ is geometrically integral, which is (iii).

Assume (iii), let $L/F$ be a field extension, and put $D=S\otimes_F L$, a domain.
If $a,b\in D$ satisfy $a^{kn}=b^n$, then $(a^k-b)\,F_{k,n}(a,b)=a^{kn}-b^n=0$,
so $a^k-b=0$ or $F_{k,n}(a,b)=0$, because $D$ is a domain.
Thus the zero ideal of $D$ is $(kn,n)$-mpdf, that is, $D$ satisfies $\MP{k}{n}$.
\end{proof}

\begin{example}\label{ex:base-change}
For a base-change failure, take $S=\mathbb R\times\mathbb R$ and
$n=3$.  Cubing is injective on $S$, so it satisfies
$\MP{k}{3}$ for every $k$.  But
$S\otimes_{\mathbb R}\mathbb C\cong\mathbb C\times\mathbb C$
fails: with $\zeta$ a primitive cube root, use $a=(1,1)$ and
$b=(\zeta,1)$.  The first factor is $(1-\zeta,0)$ and the cofactor
has second coordinate $3$.
\end{example}

\subsection{Pullbacks of domains}

The spectral theorem also gives a complete criterion for a pullback with
two domain components. The same-domain duplication case follows as a corollary.

\begin{theorem}\label{thm:domain-pullback}
Let $A,B$ be domains, let $f:A\twoheadrightarrow C$ and
$g:B\twoheadrightarrow C$ be surjective unital homomorphisms, and
suppose both kernels are nonzero. Fix $n\ge2$ and $k\ge1$, and let
$P=A\times_C B$.
Then the zero ideal of $P$ is $(kn,n)$-mpdf if and only if
\[
 \bigl(\mu_n(\operatorname{Frac}A)=\{1\}\ \text{or}\ n=0\text{ in }B\bigr)
 \quad\text{and}\quad
 \bigl(\mu_n(\operatorname{Frac}B)=\{1\}\ \text{or}\ n=0\text{ in }A\bigr).
\]
\end{theorem}
\begin{proof}
The ring $P$ is reduced as a subring of $A\times B$.  The kernels
$P_A=\ker(P\to A)$ and $P_B=\ker(P\to B)$ are prime, have product
zero, and are incomparable because both original kernels are
nonzero.  Every prime contains one of them, so they are precisely
the minimal primes.  The quotient domains are $A$ and $B$, and
Theorem~\ref{thm:spectral} gives exactly the two stated conditions.
\end{proof}

If one original kernel is zero, the pullback is isomorphic to the
other domain, so its zero ideal is automatically mpdf; this case
must be separated from the displayed two-minimal-prime criterion.

The duplication of a domain is the pullback of two copies of its
quotient map. This application is included here with the general
two-component classification.

\begin{corollary}\label{cor:domain-duplication}
Let $A$ be an integral domain, let $0\ne J\subseteq A$ be an ideal,
and let $K=\operatorname{Frac}(A)$. For $n\ge2$ and $k\ge1$, the zero
ideal of $A\bowtie J=\{(a,a+j):a\in A,\ j\in J\}$ is
$(kn,n)$-mpdf if and only if $n=0$ in $A$ or $\mu_n(K)=\{1\}$.
\end{corollary}
\begin{proof}
Use $A\bowtie J\cong A\times_{A/J}A$ in
Theorem~\ref{thm:domain-pullback}; both kernels are $J\ne0$.
\end{proof}

\begin{example}
The reduced ring $\mathbb Z\bowtie2\mathbb Z$ satisfies
$\MP{k}{3}$ for every $k\ge1$, since $\mu_3(\mathbb Q)=\{1\}$.
It never satisfies $\MP{k}{2}$, since $-1\ne1$ in $\mathbb Q$.
\end{example}

\begin{example}\label{ex:domain-pullback}
For a mixed-characteristic example, let $\zeta\in\mathbb C$ be a
primitive cube root of unity and put
\[
 P=\mathbb Z[\zeta]\times_{\mathbb F_3}\mathbb F_3[t],
\]
where the maps send $\zeta$ to $1$ and $t$ to $0$.  Both are
surjective and have nonzero kernels.  The fraction field
$\mathbb Q(\zeta)$ contains nontrivial cube roots of unity, while
$\mathbb F_3(t)$ has none, because $X^3-1=(X-1)^3$ in
characteristic three.  Also $3=0$ in the second domain.
Thus $P$ is reduced, not a domain, and satisfies
$\MP{k}{3}$ for every $k\ge1$.
\end{example}

\section{The characteristic-\texorpdfstring{$p$}{p} reduction}\label{sec:charp}

The exponent's prime-to-characteristic part has no effect on the mpdf
condition once the characteristic divides the exponent.
For $n=3$ and $k=1$, the equivalence with
$c^3=0\Rightarrow c^2=0$ follows from
\cite[Proposition~2.4 and Theorem~2.6]{Far}. For quartic pair
absorption in characteristic two, the analogous equivalence with
$c^4=0\Rightarrow c^3=0$ is \cite[Proposition~2.18]{TAHA}.
The next theorem identifies the exact characteristic part for
arbitrary exponents and mixed multipliers.

\begin{theorem}\label{thm:exact-p-part}
Let $p$ be a prime and let $S$ be a nonzero commutative ring of characteristic
$p$. Let
$n\geq2$ satisfy $p\mid n$, write $n=dq$ with
$d=p^{v_p(n)}$ and $p\nmid q$, and let $k\geq1$.
The following are equivalent:
\begin{enumerate}
\item $(0)$ is $(kn,n)$-mpdf in $S$;
\item for every $c\in S$, $c^d=0$ implies $c^{d-1}=0$;
\item $(0)$ is $(\ell d,d)$-mpdf in $S$ for every $\ell\geq1$;
\item $(0)$ is $(\ell n,n)$-mpdf in $S$ for every $\ell\geq1$.
\end{enumerate}
Whenever these conditions hold, the cofactor vanishes in every mixed-power
identity: $a^{\ell n}=b^n$ implies $F_{\ell,n}(a,b)=0$.
\end{theorem}
\begin{proof}
In characteristic $p$, $F_d(X,Y)=(X-Y)^{d-1}$, and
\[
 F_n(X,Y)=F_d(X,Y)F_q(X^d,Y^d)
         =F_q(X,Y)F_d(X^q,Y^q).
\]
Assume (i) and let $c^d=0$. Set $a=1$, $b=1+c$.
Frobenius gives $b^d=1$, hence $a^{kn}=b^n=1$.
The two factors are $-c$ and
\[
 F_{k,n}(1,1+c)=F_n(1,1+c)=q(-c)^{d-1}.
\]
The mpdf property therefore implies $c=0$ or $q(-c)^{d-1}=0$.
The integer $q$ represents a unit in $S$, so $c^{d-1}=0$ in both cases.
This proves (i)$\Rightarrow$(ii).

Assume (ii), fix $\ell\geq1$, and suppose $a^{\ell n}=b^n$.
Put $x=a^\ell$, $y=b$, and $w=x^q-y^q$. Then $w^d=0$, so $w^{d-1}=0$.
Consequently
\[
 F_{\ell,n}(a,b)=F_q(x,y)(x^q-y^q)^{d-1}=0.
\]
This proves (ii)$\Rightarrow$(iv) and the final assertion.
The same argument with $q=1$ gives (ii)$\Rightarrow$(iii).
Conversely, (iii) implies (ii) by taking $a=0$, $b=c$ in its defining
implication; and (iv) implies (i) immediately.
\end{proof}

The power-closure condition is local, although the unrestricted mpdf
condition need not be local. This is consistent with the localization
results for closed ideals in \cite[Corollary~2.9]{ABclosed}.

\begin{corollary}\label{cor:charp-local}
Let $p$ be a prime, let $S$ be a nonzero ring of characteristic $p$, and let $n\ge2$ be
divisible by $p$. For any $k\ge1$, the following are equivalent:
$S$ satisfies $\MP{k}{n}$; every $S_{\mathfrak m}$, with
$\mathfrak m\in\operatorname{Max}S$, satisfies $\MP{k}{n}$; every
$S_P$, with $P\in\Spec S$, satisfies $\MP{k}{n}$.
\end{corollary}
\begin{proof}
Localization preserves the property. Conversely, put
$d=p^{v_p(n)}$ and suppose every maximal localization has the property.
If $c^d=0$ in $S$, Theorem~\ref{thm:exact-p-part} gives
$c^{d-1}/1=0$ in every maximal localization. An element that vanishes
in all maximal localizations is zero, so the same theorem applies to $S$.
\end{proof}

Theorem~\ref{thm:exact-p-part} also applies to an ideal $I$ whenever
$p\in I$, by passing to $R/I$. It strengthens the prime-power theorem
and the one-way divisibility implication in \cite{TS} to an equivalence for all exponents with the same $p$-part.

The power-closure condition places a sharp bound on the order of every
individual nilpotent, without asserting nilpotence of the whole nilradical.

\begin{proposition}\label{prop:nilpotent-index}
Let $d\geq2$ and suppose that a commutative ring $S$ satisfies
$c^d=0\Rightarrow c^{d-1}=0$ for every $c\in S$.
If $z$ has nilpotency index $h$, meaning $z^h=0$ and $z^{h-1}\neq0$
for $z\neq0$ (with index $1$ for $z=0$), then
\[
 h\leq(d-1)\left\lceil\frac hd\right\rceil,
 \qquad\text{and consequently}\qquad h\leq(d-1)^2.
\]
Equivalently, $h=dq+r$ with integers $q\geq0$, $1\leq r\leq d-1$,
and $q+r\leq d-1$.
\end{proposition}
\begin{proof}
Put $s=\lceil h/d\rceil$. Then $(z^s)^d=0$, so
$(z^s)^{d-1}=0$, giving $h\leq(d-1)s$.
Since $s\leq(h+d-1)/d$, this yields
$dh\leq(d-1)h+(d-1)^2$, hence $h\leq(d-1)^2$.
If $d\mid h$ and $h>0$, the first inequality is impossible.
Otherwise write $h=dq+r$ with $1\leq r\leq d-1$ and use
$\lceil h/d\rceil=q+1$ to obtain $q+r\leq d-1$.
\end{proof}

The next numerical criterion is the $(d,d-1)$ specialization of the
prime-power arithmetic of closed ideals in \cite[Theorem~3.1]{ABclosed}.
We give its short valuation proof to make the characteristic-$p$
applications explicit.

\begin{proposition}\label{prop:truncated-nilpotent}
Let $K$ be a field, let $d\geq2$, and let $h\geq1$.
The ring $K[T]/(T^h)$ satisfies $c^d=0\Rightarrow c^{d-1}=0$ if and only if
$h\leq(d-1)\lceil h/d\rceil$.
In particular, the bound $(d-1)^2$ in
Proposition~\ref{prop:nilpotent-index} is attained.
If $\operatorname{char}K=p$, $d=p^e$ with $e\geq1$, and
$h=(d-1)^2$, then $K[T]/(T^h)$ is $(kn,n)$-mpdf for every $k\geq1$
and every $n$ with $v_p(n)=e$.
\end{proposition}
\begin{proof}
The necessity follows by applying the preceding proposition to the class
of $T$ (and is immediate for $h=1$). Every nonzero nilpotent in
$K[T]/(T^h)$ has the form $T^s u$ with $1\leq s<h$ and $u$ a unit.
Its $d$th power vanishes precisely when $ds\geq h$.
Thus the condition in the statement holds precisely when
$(d-1)s\geq h$ for every integer $s\geq\lceil h/d\rceil$; it suffices to
test the least such integer. For $h=(d-1)^2$,
$\lceil h/d\rceil=d-1$, giving equality.
The final assertion follows from Theorem~\ref{thm:exact-p-part}.
\end{proof}

\begin{remark}\label{quartic:characteristic-two}
If $2\in I$, the quotient form of \cite[Proposition~2.18]{TAHA}
says that $I$ is 4-pdf-absorbing exactly when
$c^4\in I\Rightarrow c^3\in I$ for every $c\in R$.
Theorem~\ref{thm:exact-p-part} identifies this condition with
$(4k,4)$-mpdf for some, equivalently every, $k\ge1$; the direct
comparison is also given in \cite[Corollary~3.21]{TS}.
For a field $F$ of characteristic two, the known quartic case
\cite[Corollary~2.19]{TAHA} is
\[
 (T^e)\subset F[T]\text{ is 4-pdf-absorbing}
 \quad\Longleftrightarrow\quad e\in\{1,2,3,5,6,9\}.
\]
This is the $d=4$ case of Proposition~\ref{prop:truncated-nilpotent};
the same list therefore governs $(4k,4)$-mpdf for every $k$.
\end{remark}

\begin{example}
The ring $\mathbb F_3[T]/(T^4)$ is $(3k,3)$-mpdf for every $k\geq1$,
although the image of $T$ has nilpotency index $4>3-1$.
Thus the correct uniform bound is $(d-1)^2$, not $d-1$.
The nilradical itself need not be nilpotent: for example, in
$S=\mathbb F_p[X_i:i\geq1]/(X_i^p:i\geq1)$ every nilpotent $c$ satisfies
$c^p=0$, so $S$ is $(kp^2,p^2)$-mpdf for all $k\geq1$; yet the ideal
generated by the images of all $X_i$ has nonzero products of arbitrarily
large length. Indeed, a product of distinct variables survives in this
monomial quotient, and Frobenius sends every polynomial with zero constant
term to zero.
\end{example}

\section{Square-zero extensions and idealization ideals}\label{sec:squarezero}

We first treat extensions whose reduced quotient is a domain.  The extension
need not split, and the mixed-characteristic case retains information that
is invisible in the quotient and its module alone.
For $n=2$, such a ring satisfies $\MP{k}{2}$ exactly when its
square-zero kernel is zero: the power-closure consequence makes the
ring reduced, and the quotient is already a domain. Thus the substantial
extension phenomena begin at $n=3$.

\Needspace{9\baselineskip}
\begin{theorem}\label{thm:squarezero-domain}
Let $S$ be a commutative ring, let $N$ be an ideal with $N^2=0$, and suppose
that $A=S/N$ is an integral domain. Fix integers $n\geq3$ and $k\geq1$.
Regard $N$ as an $A$-module in the natural way.
\begin{enumerate}
\item If $n\cdot1_A\neq0$, then $(0)$ is $(kn,n)$-mpdf in $S$ if and only
if $N$ is torsion-free over $A$.
\item Suppose $\operatorname{char}A=p>0$, $p\mid n$, and
$d=p^{v_p(n)}\geq3$. Then $(0)$ is $(kn,n)$-mpdf in $S$ if and only if
$n\cdot1_S=0$.
\item Suppose $\operatorname{char}A=2$ and $v_2(n)=1$. Then $(0)$ is
$(kn,n)$-mpdf in $S$ if and only if $N=\{0,2\cdot1_S\}$.
Here the latter set may consist of one element.
\end{enumerate}
In particular, in all these extensions the condition is independent of $k$.
\end{theorem}

\begin{proof}
If $v\in N$, expansion using $v^2=0$ gives
\[
 (x+v)^n=x^n+n x^{n-1}v,
 \qquad
 F_n(x+v,x)=n x^{n-1}+\binom n2 x^{n-2}v.
\]

For (i), assume first that $S$ has the stated mpdf property. Let
$0\neq r\in A$ and $v\in N$ satisfy $rv=0$, and choose a lift $a\in S$ of
$r$. Set $b=a^k+v$. Then $a^{kn}=b^n$, since
$n a^{k(n-1)}v=0$. The reduction of $F_n(a^k,b)$ is
$n r^{k(n-1)}\neq0$ in $A$. Thus the first factor must vanish, giving
$v=0$. Hence $N$ is torsion-free.

Conversely, suppose $N$ is torsion-free and $a^{kn}=b^n$. Put
$x=a^k$, $y=b$, $D=x-y$, and $C=F_n(x,y)$, so $DC=0$.
If $\overline D\neq0$, then $\overline C=0$, hence $C\in N$; the equality
$\overline D C=0$ and torsion-freeness give $C=0$.
If $\overline D=0$ but $\overline x\neq0$, then
$\overline C=n\overline x^{n-1}\neq0$, and torsion-freeness gives $D=0$.
The only remaining case is $\overline x=\overline y=0$. Since $A$ is a
domain, $a,b\in N$, and every monomial of $F_{k,n}(a,b)$ has total degree
at least $n-1\geq2$. Hence $F_{k,n}(a,b)=0$.

For (ii), note first that $p\in N$, so $p^2=0$ in $S$ and $pN=0$.
The integer $\binom n2$ is divisible by $p$: this is immediate for odd
$p$, and for $p=2$ the hypothesis $d\geq3$ implies $4\mid n$.
If $n\cdot1_S\neq0$, set $v=n\cdot1_S\in N$ and consider $a=1$ and
$b=1+v$. Then $a^{kn}=b^n=1$, while $a^k-b=-v\neq0$ and
$F_{k,n}(a,b)=n\cdot1_S\neq0$, a contradiction.

For sufficiency write $n=dq$, with $p\nmid q$, and suppose $n\cdot1_S=0$.
If $d=p$, then $p$ is odd and $p\cdot1_S=0$: indeed $q$ is a unit in
$S$, because $p^2=0$ and $\gcd(p,q)=1$. In characteristic $p$, if
$x^n=y^n$, then $(x^q-y^q)^p=0$; reduction to the domain $A$ shows that
$x^q-y^q\in N$, and consequently
\[
 F_n(x,y)=F_q(x,y)(x^q-y^q)^{p-1}=0
\]
because $p-1\geq2$.
If $d\geq p^2$, again $w=x^q-y^q\in N$ after reducing $x^n=y^n$ modulo
$N$. Writing $u=y^q$ gives
\[
 F_d(u+w,u)=d u^{d-1}+\binom d2 u^{d-2}w=0.
\]
Here $d\cdot1_S=0$ because $p^2=0$, and $p\mid\binom d2$, so the second
term vanishes because $pN=0$. Since
$F_n(x,y)=F_q(x,y)F_d(x^q,y^q)$, the required cofactor is zero.

For (iii), write $n=2q$ with $q$ odd. Again $2\in N$, $4=0$, and $2N=0$.
For any $v\in N$ the pair $a=1$, $b=1+v$ satisfies $a^{kn}=b^n=1$,
and its two factors are $-v$ and
$n+\binom n2v=2+v$. Thus every nonzero $v\in N$ must equal $2$, proving
necessity.

Conversely, suppose $N=\{0,2\}$. If $N=0$, then $S=A$ is a domain,
so suppose $|N|=2$. Let $x=a^k$, $y=b$ satisfy $x^n=y^n$.
In $A$ one has $\overline x^q=\overline y^q$, because the square map is
injective on a domain of characteristic two. If
$\overline x\neq\overline y$, then $F_q(x,y)\in N$ and
$x^q+y^q\in N$, whence
$F_n(x,y)=F_q(x,y)(x^q+y^q)=0$.
If $\overline x=\overline y$, either $x=y$ or $x-y=2$. In the latter case
\[
 F_n(x,y)=n y^{n-1}+\binom n2 y^{n-2}\,2
         =2y^{n-2}(y+1)=0.
\]
For the last equality, multiplication by $y$ on the two-element module
$N$ is either zero or the identity, so $2y(y+1)=0$; use $n-2\geq1$.
This completes the proof.
\end{proof}

\begin{corollary}\label{cor:domain-idealization}
Let $A$ be an integral domain, let $M$ be an $A$-module, and let
$n\geq3$, $k\geq1$. Then the zero ideal of $A\ltimes M$ is $(kn,n)$-mpdf
precisely in the following cases:
\begin{enumerate}
\item $n\cdot1_A\neq0$ and $M$ is torsion-free;
\item $\operatorname{char}A=p\mid n$ and $p^{v_p(n)}\geq3$;
\item $\operatorname{char}A=2$, $v_2(n)=1$, and $M=0$.
\end{enumerate}
\end{corollary}
\begin{proof}
Apply Theorem~\ref{thm:squarezero-domain} to $N=0\ltimes M$.
The characteristic of the idealization equals that of $A$.
\end{proof}

\subsection{Ideals that need not be homogeneous}

For an $R$-module $M$, the idealization $R\ltimes M$ has multiplication
$(r,u)(s,v)=(rs,rv+su)$. We use the standard idealization framework of
\cite{AW}. For the strict cubic condition, transfer through saturated
ideals $I\ltimes M$ is established in \cite[Theorem~4.6]{Dr26};
related implications for nonzero $I$ appear in
\cite[Theorem~2.20]{Far}. The cdf-primary analogue is
\cite[Theorem~7(3)]{YSO}. These saturated ideals have quotient $R/I$.
The quotients below can retain a nonzero square-zero kernel and
depend on the defining linear map of a nonhomogeneous ideal.
The following description makes that information explicit.

\begin{lemma}\label{lem:compatible-idealization}
Let $R$ be a commutative ring and $M$ an $R$-module. Every proper ideal $L$ of
$R\ltimes M$ is uniquely specified by a proper ideal $P$ of $R$, a submodule
$N$ of $M$ with $PM\subseteq N$, and an $R$-linear map
$\phi:P\to M/N$, via
\[
 L=L(P,N,\phi)=\{(r,u):r\in P,\ u+N=\phi(r)\}.
\]
Conversely, every such triple defines an ideal. The map $\phi$ kills
$P^2$ and may equivalently be regarded as an $R/P$-linear map
$P/P^2\to M/N$. Moreover, if $S=(R\ltimes M)/L$, the natural
surjection $S\to R/P$ has a square-zero kernel canonically identified
with the $R/P$-module $M/N$. If an integer $t$ belongs to $P$, its image
$t\cdot1_S$ corresponds to $-\phi(t)$ in the kernel.
\end{lemma}
\begin{proof}
For a given $L$, take $P$ to be its first-coordinate projection and put
$N=\{u:(0,u)\in L\}$. If $(r,u)\in L$, multiplication by $(0,m)$ gives
$rm\in N$, so $PM\subseteq N$. The coset $u+N$ depends only on $r$;
define it to be $\phi(r)$. Addition and multiplication by $(a,0)$ prove
$R$-linearity. Conversely these properties show directly that the
displayed subset is an ideal. Since $P$ kills $M/N$, linearity implies
$\phi(P^2)=0$. The quotient description follows from first-coordinate
projection. If $u+N=\phi(t)$, then $(t,u)\in L$, so
$(t,0)+L=-(0,u)+L$, which gives the final assertion.
\end{proof}

\begin{corollary}\label{cor:arbitrary-idealization-mpdf}
In Lemma~\ref{lem:compatible-idealization}, suppose $P$ is prime and fix
$n\geq3$, $k\geq1$. Put $A=R/P$ and $E=M/N$.
\begin{enumerate}
\item If $n\notin P$, then $L$ is $(kn,n)$-mpdf if and only if $E$ is
torsion-free over $A$.
\item If $\operatorname{char}A=p\mid n$ and $p^{v_p(n)}\geq3$, then $L$
is $(kn,n)$-mpdf if and only if $\phi(n)=0$.
\item If $\operatorname{char}A=2$ and $v_2(n)=1$, then $L$ is
$(kn,n)$-mpdf if and only if $E=\{0,\phi(2)\}$.
\end{enumerate}
\end{corollary}
\begin{proof}
Use the quotient criterion and apply Theorem~\ref{thm:squarezero-domain}
to the square-zero quotient in Lemma~\ref{lem:compatible-idealization}.
\end{proof}

\begin{example}\label{ex:nonsplit-characteristic-four}
Let $R=\mathbb Z$, $M=\mathbb F_2$, $P=2\mathbb Z$, $N=0$, and
$\phi(2a)=\overline a$. Then $L=L(P,0,\phi)$ is not homogeneous and
$(R\ltimes M)/L\cong\mathbb Z/4\mathbb Z$, via
$(r,u)\mapsto r+2u\pmod4$. Hence $L$ is $(6k,6)$-mpdf for every $k\geq1$.
In contrast, the homogeneous ideal $L_0=2\mathbb Z\ltimes0$ has quotient
$\mathbb F_2\ltimes\mathbb F_2$, which is not $(6k,6)$-mpdf for any $k$.
Indeed, in the latter quotient take $a=(1,0)$ and $b=(1,1)$: their
required powers agree, while both factors equal $(0,1)$.
Thus identical $P$ and $N$ need not give identical mpdf behavior; the
extension map $\phi$ matters.
\end{example}

\subsection{Homogeneous quartic ideals over primary bases}
\label{quartic:primary-idealization-section}

The prime-projection classification above retains a square-zero
kernel. A homogeneous idealization can also have a primary base
ideal whose radical has nonzero square. Under regularity of the
module action, \cite[Theorem~5.5]{TAHA} characterizes the quartic
pair-product condition for such ideals when $2,3\notin P$.
The primary quartic comparison in \cite[Theorem~3.17]{TS} states that, for a
$P$-primary ideal with $2\notin P$, the 4-pdf and $(4,4)$-mpdf
conditions are both equivalent to $P^3\subseteq I$.
It extends \cite[Theorem~4.6]{TAHA} by removing the restriction on
$3$. Applying that comparison gives the following version of the
homogeneous criterion.

\begin{corollary}\label{quartic:homogeneous-idealization}
Let $I$ be a $P$-primary ideal of $R$ with $2\notin P$. Let $M$
be an $R$-module and $N\subseteq M$ a submodule with $IM\subseteq N$.
Suppose every element of $R\setminus P$ acts injectively on $M/N$.
For $J=I\ltimes N\subset R\ltimes M$, the following are equivalent:
\begin{enumerate}
\item $J$ is 4-pdf-absorbing.
\item $J$ is $(4,4)$-mpdf.
\item $P^3\subseteq I$ and $P^2M\subseteq N$.
\end{enumerate}
\end{corollary}
\begin{proof}
Put $A=R/I$, $\mathfrak p=P/I$, and $E=M/N$, an $A$-module.
Then
\[
 B=A\ltimes E\cong(R\ltimes M)/J.
\]
Its nilradical is $K=\mathfrak p\ltimes E$, with
$B/K\cong R/P$ a domain. The zero ideal of $B$ is $K$-primary:
if $(a,u)\notin K$ and $(a,u)(b,v)=0$, then $ab=0$ in $A$.
Primaryness of $I$ gives $b=0$, since $a\notin\mathfrak p$.
The second coordinate then gives $av=0$, so $v=0$ by the
module hypothesis. Thus every element outside $K$ is regular.

Also $2\notin K$. By the primary quartic comparison in
\cite[Theorem~3.17]{TS},
the zero ideal of $B$ satisfies either of the two quartic
conditions exactly when $K^3=0$. Multiplication gives
\[
 K^3=\mathfrak p^3\ltimes\mathfrak p^2E.
\]
Indeed, products of three elements give one containment; products
of $(r,0),(s,0),(t,0)$ and of $(r,0),(s,0),(0,e)$ generate the
two indicated coordinates and give the reverse containment.
Consequently $K^3=0$ is equivalent to
$P^3\subseteq I$ and $P^2M\subseteq N$.
Testing the two conditions in the quotient $B$ proves the result.
\end{proof}

Taking $M=R$ and $N=I$ gives the useful specialization
\[
 I\ltimes I\text{ satisfies either quartic condition}
 \quad\Longleftrightarrow\quad P^2\subseteq I.
\]
Indeed, every element of $R\setminus P$ acts injectively on $R/I$,
and $P^2\subseteq I$ implies $P^3\subseteq I$.
This extends \cite[Corollary~5.6]{TAHA} under the sole restriction
$2\notin P$.

The quotient-primary argument and ideal-power calculation are those
used in \cite[Theorem~5.5]{TAHA}; the corollary removes its
restriction on $3$ and identifies the common strict quartic case.
Here $I$ need not be prime, and $K^2$ can be nonzero. In
Corollary~\ref{cor:arbitrary-idealization-mpdf}, the first-coordinate
projection is prime, but the ideal may be nonhomogeneous and
the mixed parameters range over $n\ge3$ and $k\ge1$.
Theorem~\ref{thm:squarezero-domain} also treats nonsplit extensions.
These classifications therefore share a quartic case while
retaining different hypotheses on the base and extension.

\subsection{Geometric stability of square-zero extensions}

The reducedness hypothesis in Theorem~\ref{thm:geometric} is essential for
$n\ge3$. The next result identifies exactly which square-zero thickenings
of geometrically reduced algebras retain the mpdf property after every
extension of the ground field. The extension need not split.

\begin{lemma}\label{lem:geometric-torsionfree}
Let $F$ be a field, let $A$ be a geometrically integral $F$-algebra,
and let $M$ be a torsion-free $A$-module. For every field extension
$L/F$, the $A\otimes_F L$-module $M\otimes_F L$ is torsion-free.
\end{lemma}
\begin{proof}
Put $A_L=A\otimes_F L$, which is a domain. Every finitely generated
$A$-submodule $M_0$ of $M$ embeds in a free $A$-module of finite
rank, by the standard denominator-clearing argument
\cite[Tag 0AUU]{Stacks}: embed
$M_0$ in $M_0\otimes_A\operatorname{Frac}(A)$, choose a basis of this
finite-dimensional vector space, and clear the denominators of a finite
set of generators. Tensoring the resulting injection $M_0\hookrightarrow
A^r$ over $F$ with the flat $F$-module $L$ gives an injection
$M_0\otimes_F L\hookrightarrow A_L^r$. Thus $M_0\otimes_F L$ is
torsion-free over $A_L$. Since $M$ is the directed union of its finitely
generated submodules, and tensoring over $F$ preserves their injections,
$M\otimes_F L$ is a directed union of such torsion-free $A_L$-modules.
It is therefore torsion-free.
\end{proof}

\begin{theorem}\label{thm:geometric-squarezero}
Let $F$ be a field, let $S$ be a nonzero commutative $F$-algebra, and
let $N$ be an ideal of $S$ with $N^2=0$. Suppose that $A=S/N$ is
geometrically reduced over $F$. Fix $n\ge3$ and $k\ge1$, and suppose
$\operatorname{char}F\nmid n$. The following are equivalent:
\begin{enumerate}
\item $S\otimes_F L$ satisfies $\MP{k}{n}$ for every field extension
$L/F$;
\item $S\otimes_F\overline F$ satisfies $\MP{k}{n}$;
\item $A$ is geometrically integral over $F$, and $N$ is torsion-free
as an $A$-module.
\end{enumerate}
When these conditions hold, $S\otimes_F L$ satisfies
$\MP{\ell}{n}$ for every field extension $L/F$ and every $\ell\ge1$.
\end{theorem}
\begin{proof}
The implication (i)$\Rightarrow$(ii) is immediate. Assume (ii), and
write $S_{\overline F}=S\otimes_F\overline F$ and
$N_{\overline F}=N\otimes_F\overline F$. Flatness over $F$ gives
$S_{\overline F}/N_{\overline F}\cong A\otimes_F\overline F$.
This quotient is reduced and $N_{\overline F}^2=0$, so
$N_{\overline F}$ is the nilradical of $S_{\overline F}$.
Corollary~\ref{thm:radical-descent} shows that
$A\otimes_F\overline F$ satisfies $\MP{k}{n}$.
Theorem~\ref{thm:geometric} now implies that $A$ is geometrically
integral. The injection $S\hookrightarrow S_{\overline F}$ and
Lemma~\ref{lem:transfers}(ii) imply that $S$ satisfies $\MP{k}{n}$.
Since $A$ is a domain and $n\cdot1_A\ne0$, Theorem~\ref{thm:squarezero-domain}(i)
then shows that $N$ is torsion-free over $A$. This proves (iii).

Assume (iii), and let $L/F$ be any field extension. Tensoring
$0\to N\to S\to A\to0$ over $F$ with $L$ gives a square-zero
extension of the domain $A_L=A\otimes_F L$ by
$N_L=N\otimes_F L$. Lemma~\ref{lem:geometric-torsionfree} makes
$N_L$ torsion-free over $A_L$. Since $n\cdot1_{A_L}\ne0$,
Theorem~\ref{thm:squarezero-domain}(i) gives $\MP{\ell}{n}$ for
$S\otimes_F L$ for every $\ell\ge1$. This proves (i) and the final
assertion.
\end{proof}

\begin{example}
Let $F$ be a field and let $n\ge3$ be prime to its characteristic.
The nonreduced algebra $F[\varepsilon]/(\varepsilon^2)$ satisfies
$\MP{k}{n}$ after every extension of $F$, for every $k\ge1$:
its quotient by $(\varepsilon)$ is $F$, and $(\varepsilon)$ is a
one-dimensional vector space over $F$. Thus universal mpdf stability
does not by itself imply geometric integrality once $n\ge3$.
\end{example}

\section{Valuation and Dedekind-domain classifications}\label{sec:valuation}

A discrete valuation ring (DVR) is a local principal ideal domain that
is not a field. We write $v$ for its normalized valuation and $\pi$
for a uniformizer, so $v(\pi)=1$ and $v(0)=\infty$.
The term \emph{tame} below means that the integer $n$ is a unit in the
DVR. The following criterion separates the usual power-closure
restriction from the additional restriction imposed by the mixed power.

\begin{theorem} \label{thm:tame-dvr}
Let $V$ be a discrete valuation ring with uniformizer $\pi$, let $n\ge2$
be a unit of $V$, and let $k,e\ge1$. Put $h=\lceil e/n\rceil$.
Then $\pi^eV$ is $(kn,n)$-mpdf if and only if
\[
 h\le k\qquad\text{and}\qquad e\le(n-1)h.
\]
Equivalently, the admissible exponents are precisely
\[
 \bigcup_{h=1}^{\min\{k,n-1\}}
 \{n(h-1)+1,\ldots,(n-1)h\}.
\]
\end{theorem}
\begin{proof}
Write $v$ for the normalized valuation, with $v(0)=\infty$.
For necessity, set $a=0$ and $b=\pi^h$. Then $b^n\in\pi^eV$.
The conclusion forces $b^{n-1}\in\pi^eV$: the other alternative,
$b\in\pi^eV$, also implies this. Hence $e\le(n-1)h$.
Suppose next that $e>nk$, and put $s=e-nk\ge1$,
$a=\pi$, and $b=\pi^k(1+\pi^s)$. Because $n$ is a unit,
$F_n(1,1+\pi^s)$ is a unit. Consequently
\[
 v(a^k-b)=k+s=e-(n-1)k<e,
 \qquad v(F_{k,n}(a,b))=(n-1)k<e,
\]
whereas their product has valuation $nk+s=e$. This contradicts the
mpdf condition. Therefore $e\le nk$, equivalently $h\le k$.

Conversely, assume both inequalities and let $a^{kn}-b^n\in\pi^eV$.
If $a$ is a unit, then so is $b$. Put $x=a^k$. If $x-b$ is a unit,
the product identity gives $F_n(x,b)\in\pi^eV$. If $x-b$ is not a
unit, then $x\equiv b\pmod\pi$, so
$F_n(x,b)\equiv nb^{n-1}\not\equiv0\pmod\pi$. Thus the cofactor is a
unit, and $x-b\in\pi^eV$.
If $a$ is a nonunit, then $a^{kn}\in\pi^eV$, since $nk\ge e$.
It follows that $b^n\in\pi^eV$, so $v(b)\ge h$. Also $v(a^k)\ge k\ge h$.
Every summand of $F_n(a^k,b)$ has valuation at least $(n-1)h\ge e$.
This proves sufficiency.
Finally, $h=\lceil e/n\rceil$ means $n(h-1)+1\le e\le nh$.
Together with $e\le(n-1)h$, this interval is nonempty precisely when
$h\le n-1$, giving the displayed union.
\end{proof}

The valuation criterion realizes every possible threshold before
stabilization, rather than only the first separation between $k=1$ and
$k=2$.

\begin{corollary} \label{cor:hierarchy}
Under the hypotheses of Theorem~\ref{thm:tame-dvr}, put
$d=\min\{k,n-1\}$. Exactly $dn-d(d+1)/2$ positive exponents $e$ give
$(kn,n)$-mpdf ideals $\pi^eV$. For fixed $k$, the largest exponent is $(n-1)d$. The list stabilizes at
$k=n-1$, when its largest exponent is $(n-1)^2$. For each $r\in\{1,\ldots,n-1\}$,
$\pi^{r(n-1)}V$ is $(kn,n)$-mpdf precisely when $k\ge r$.
\end{corollary}
\begin{proof}
The interval indexed by $h$ contains $n-h$ integers, and the intervals are
disjoint, so their total cardinality is $\sum_{h=1}^d(n-h)$.
For $1\le r\le n-1$, we have
$\lceil r(n-1)/n\rceil=r$, and $r(n-1)=(n-1)r$.
Apply the theorem.
\end{proof}

\begin{remark}
The case $k=1$ says that the admissible exponents are $1,\ldots,n-1$.
For prime $n$, this recovers the tame characteristic-zero PID result of
\cite[Theorem~3]{Yes}, after localization. The general theorem allows
arbitrary $n$, arbitrary $k$, and arbitrary residue field, and determines
the exact dependence on $k$.
\end{remark}

\Needspace{16\baselineskip}
\begin{remark}\label{quartic:dvr-comparison}
For a DVR $V$ whose residue characteristic is different from two,
the primary quartic comparison in \cite[Theorem~3.17]{TS} makes $\pi^eV$
4-pdf-absorbing exactly when $e\le3$. This uses only $2\notin\pi V$,
so it includes residue characteristic three, beyond the hypotheses
of \cite[Theorem~4.6]{TAHA}. Theorem~\ref{thm:tame-dvr}, with $n=4$,
gives the exact comparison:
\begin{center}
\begin{tabular}{@{}ll@{}}
\toprule
Property of $\pi^eV$ & Permitted positive exponents $e$\\
\midrule
4-pdf-absorbing & $1,2,3$\\
$(4,4)$-mpdf & $1,2,3$\\
$(8,4)$-mpdf & $1,2,3,5,6$\\
$(4k,4)$-mpdf, $k\ge3$ & $1,2,3,5,6,9$\\
\bottomrule
\end{tabular}
\end{center}
The extra levels in the last two rows are permitted by the mixed
multiplier; they do not satisfy the unmixed pair-product condition.
\end{remark}

\begin{example}\label{quartic:mixed-separation}
Let $I=5^5\Z$. In the DVR $V=\Z_{(5)}$, take $n=4$, $e=5$,
and $k\ge2$. Since $h=\lceil5/4\rceil=2\le k$ and $5\le3h$,
Theorem~\ref{thm:tame-dvr} makes $5^5V$ $(4k,4)$-mpdf.
Its contraction is $5^5V\cap\Z=I$. The induced embedding
$\Z/I\hookrightarrow V/5^5V$ and the quotient and subring
properties in Lemma~\ref{lem:transfers} show that $I$ is
$(4k,4)$-mpdf for every $k\ge2$.

Nevertheless $I$ is not 4-pdf-absorbing. For $a=5$, $b=10$,
\[
 \begin{aligned}
 a^4-b^4&=-9375=-3\cdot5^5\in I,\\
 a^2-b^2&=-75,\\
 (a^2+b^2)(a-b)&=-625,\\
 (a^2+b^2)(a+b)&=1875.
 \end{aligned}
\]
None of the last three integers belongs to $I$.
Conversely, $(0)\subset\Q\times\Q$ is 4-pdf-absorbing but is
not $(4k,4)$-mpdf for any $k$, as in
Example~\ref{quartic:radical-examples}. Thus, for each $k\ge2$,
the mixed strict and unmixed pair-product conditions are incomparable.
\end{example}

When the exponent itself is prime, the residue characteristic may equal
that exponent. The next criterion includes arbitrary ramification.

\Needspace{10\baselineskip}
\begin{theorem} \label{thm:ramified-prime}
Let $V$ be a discrete valuation ring of characteristic zero with
uniformizer $\pi$ and residue characteristic $p>0$, and put $E=v_\pi(p)$.
For $e,k\ge1$, the following are equivalent:
\begin{enumerate}
\item $\pi^eV$ is $(kp,p)$-mpdf;
\item $e\le E$ and $e\le(p-1)\lceil e/p\rceil$.
\end{enumerate}
In particular the property is independent of $k$ in this ramified case.
\end{theorem}
\begin{proof}
Write $h=\lceil e/p\rceil$. The choice $a=0$, $b=\pi^h$ gives the
necessary inequality $e\le(p-1)h$, exactly as above.
Suppose this holds but $e>E$. Set $s=\max\{h,e-E\}$, and take
$a=1$, $b=1+\pi^s$. Since $e\ge2$, both $h$ and $e-E$ are strictly less
than $e$, so $s<e$. The polynomial identity
\[
 F_p(1,1+Z)=\sum_{j=1}^p\binom pj Z^{j-1}
\]
shows that $F_p(1,1+\pi^s)$ has valuation exactly $E$.
Indeed, its constant term is $p$, every middle term has valuation
at least $E+s>E$, and the final term has valuation
$(p-1)s\ge(p-1)h\ge e>E$.
The linear factor has valuation $s<e$, while the product has valuation
$s+E\ge e$. Neither factor belongs to $\pi^eV$, a contradiction.
Thus $e\le E$.

Conversely, suppose $e\le E$ and $e\le(p-1)h$.
The quotient $V/\pi^eV$ has characteristic $p$. Modulo $\pi^eV$,
\[
 a^{kp}-b^p=(a^k-b)^p,
 \qquad F_{k,p}(a,b)=(a^k-b)^{p-1}.
\]
If the first expression is zero, then $v(a^k-b)\ge h$, so the
second belongs to $\pi^eV$.
\end{proof}

\begin{remark}
Ye\c{s}ilot \cite[Remark~5]{Yes} asks for the classification of $(\pi^e)$ when $\pi$ divides
$p$ in a characteristic-zero PID.
Theorem~\ref{thm:ramified-prime} answers that question. Indeed localization
at $(\pi)$ does not alter the quotient modulo $\pi^e$, and $E$ is the
ramification valuation of $p$. The answer requires no finiteness
hypothesis on the residue field.
\end{remark}

\subsection{Unramified mixed-characteristic discrete valuations}

For arbitrary composite exponents, a complete criterion is available when
$v(p)=1$. We first prove the binomial valuation formula needed to
construct unit counterexamples.

\begin{lemma}\label{lem:principal-unit}
Let $V$ be a discrete valuation ring of characteristic zero and residue
characteristic $p>0$, with normalized valuation $v$ satisfying $v(p)=1$.
Let $\ell\ge1$ and $u\in V^\times$, and assume either that $p$ is
odd or that $\ell\ge2$. For $n\ge1$ and $t=v_p(n)$, one has
\[
 v\bigl((1+p^\ell u)^n-1\bigr)=\ell+t.
\]
\end{lemma}
\begin{proof}
If $p\nmid n$, then $n$ is a unit in $V$. The linear term in the
binomial expansion has valuation $\ell$, whereas every higher term
has valuation at least $2\ell>\ell$, proving the formula.
For $z=1+p^ru'$ with $u'\in V^\times$, where $r\ge1$ for odd $p$
and $r\ge2$ for $p=2$, the first term of $z^p-1$ has valuation $r+1$.
For odd $p$, each middle term has valuation at least $1+2r>r+1$,
and the last term has valuation $pr>r+1$. For $p=2$, the last term
has valuation $2r>r+1$. Thus $v(z^p-1)=r+1$.
Write $n=p^tq$ with $p\nmid q$. Apply the first calculation to the
$q$th power, then the second calculation successively $t$ times.
\end{proof}

The local integer classification extends to every unramified DVR of
mixed characteristic. Its characteristic-four exception depends on the
residue field and is therefore not a purely numerical ramification
phenomenon.

\Needspace{10\baselineskip}
\begin{theorem}\label{thm:unramified-dvr}
Let $V$ be a discrete valuation ring of characteristic zero and residue
characteristic $p>0$, with normalized valuation $v$ satisfying $v(p)=1$.
Let $K=V/pV$, let $n\ge2$ and $k,e\ge1$, and put $t=v_p(n)$.
Then $p^eV$ is $(kn,n)$-mpdf precisely in the following cases:
\begin{enumerate}
\item $t=0$, $h=\lceil e/n\rceil\le k$, and $e\le(n-1)h$;
\item $t>0$ and $e\le t$;
\item $p=2$, $t=1$, $e=2$, $n\ge6$, and $K\cong\mathbb F_2$.
\end{enumerate}
\end{theorem}
\begin{proof}
The tame case is Theorem~\ref{thm:tame-dvr}. Suppose $t>0$, and write
$n=p^tq$ with $p\nmid q$. If $e\le t$, let $x=a^k$, $y=b$ satisfy
$x^n-y^n\in p^eV$. When $x$ is a nonunit, so is $y$, and every term
of $F_n(x,y)$ lies in $p^{n-1}V\subseteq p^eV$, because $n-1\ge t$.
When $x$ and $y$ are units and $x-y$ is a unit, the product identity
forces $F_n(x,y)\in p^eV$. In the remaining case $x\equiv y\pmod p$,
and the factorization
\[
 F_n(x,y)=F_q(x,y)\prod_{j=0}^{t-1}
 F_p(x^{qp^j},y^{qp^j})
\]
shows that $F_n(x,y)\in p^tV\subseteq p^eV$, since every factor in
the product is in $pV$. This proves sufficiency in (ii).

For odd $p$ and $e>t$, the pair $a=1$, $b=1+p^{e-t}$ gives factors
of valuations $e-t$ and $t$, with their product of valuation $e$,
by Lemma~\ref{lem:principal-unit}. Integer valuations agree with the
normalized valuation of $V$, since $v(p)=1$ and integers prime to $p$
are units in $V$. Thus this pair is a counterexample. The same argument
works for $p=2$ when $e\ge t+2$. If $p=2$, $e=t+1$, and $t\ge2$,
the pair $a=1$, $b=5$ gives factor valuations $2$ and $t$, both less
than $e$, with product valuation $t+2\ge e$.

Only $p=2$, $t=1$, $e=2$ remains. For $n=2$, the pair $a=0$, $b=2$
violates the condition modulo $4V$. For $n\ge6$, apply
Theorem~\ref{thm:squarezero-domain}(iii) to $S=V/4V$ and
$N=2V/4V$. Multiplication by $2$ identifies $K=V/2V$ with $N$ as
$K$-modules. Since $2\ne0$ in $S$, the condition
$N=\{0,2\cdot1_S\}$ holds exactly when $K\cong\mathbb F_2$.
This is (iii), and it completes the classification.
\end{proof}

\begin{example}\label{ex:unramified-residue-obstruction}
The characteristic-four exception is sensitive to the residue field.
In the setting of Theorem~\ref{thm:unramified-dvr}, suppose $p=2$ and
$K\cong\mathbb F_4$. Choose $u\in V$ whose residue is neither $0$ nor
$1$. In $S=V/4V$, let $a=1$ and $b=1+2u$. For every $k\ge1$,
$a^{6k}=b^6=1$, but
\[
 a^k-b=-2u\ne0,\qquad F_{k,6}(a,b)=2+2u\ne0.
\]
The two nonvanishing assertions follow from the chosen residue of $u$.
Thus $4V$ fails the $(6k,6)$-mpdf condition, although $4\mathbb Z$
satisfies it for every $k$.
\end{example}

\subsection{A nondiscrete valuation contrast}

The discrete exponent ranges have no analogue for a nondiscrete rank-one
valuation. The following consequence of the necessary power-closure
condition gives the precise contrast; it uses the closed-ideal viewpoint
of \cite{ABclosed}.

\begin{corollary}\label{cor:nondiscrete-valuation}
Let $V$ be a valuation domain with surjective valuation
$v:\operatorname{Frac}(V)^\times\to\Gamma$, where $\Gamma$ is a dense
subgroup of $(\mathbb R,+)$, and let $M$ be its maximal ideal.
For every $n\ge2$ and $k\ge1$, the proper $(kn,n)$-mpdf ideals of $V$
are exactly $(0)$ and $M$.
\end{corollary}
\begin{proof}
Both listed ideals are prime, hence mpdf. Suppose $I$ is a nonzero
proper ideal different from $M$. Choose $z\in M\setminus I$.
If $y\in I\setminus\{0\}$, then $v(y)>v(z)$; otherwise
$z/y\in V$ would imply $z\in I$. Thus
\[
 \alpha=\inf\{v(y):y\in I\setminus\{0\}\}
\]
is finite and strictly positive. Every element of $V$ having value
strictly greater than $\alpha$ belongs to $I$: choose an element of
$I$ of smaller value and divide. Every element of value strictly less
than $\alpha$ lies outside $I$. Density gives $\gamma\in\Gamma$ with
$\alpha/n<\gamma<\alpha/(n-1)$; choose $c\in V$ with $v(c)=\gamma$.
Then $c^n\in I$ but $c^{n-1}\notin I$, contradicting
Lemma~\ref{lem:transfers}(iii).
\end{proof}

\subsection{Dedekind domains and number fields}

Here a Dedekind domain means a Noetherian integrally closed domain
of Krull dimension one. Every nonzero proper ideal factors uniquely
into powers of distinct maximal ideals, and the localizations of the
domain at those maximal ideals are DVRs; see \cite[Chapter~9]{AM}.
The product criterion retains the information that a purely local
classification would lose.
The finite cubic compatibility mechanism appears in
\cite[Propositions~14 and~15]{YSO} for cdf-primary ideals.
We use the strict mixed-power criterion of \cite{TS}, which allows
arbitrary exponents and arbitrary index sets. The companion paper
also explains the finite cubic overlap for nonradical strict factors.

\begin{lemma}\label{lem:product}
Let $n\ge2$ and $k\ge1$, let $\{R_i\}_{i\in\Lambda}$ be a nonempty family of rings, and let
$I_i$ be a proper ideal of $R_i$ for each $i$. Then $\prod_iI_i$ is $(kn,n)$-mpdf in
$\prod_iR_i$ if and only if the following conditions hold.
\begin{enumerate}
\item Each $I_i$ is $(kn,n)$-mpdf in $R_i$.
\item Whenever $a,b\in R_i$ satisfy $a^{kn}-b^n\in I_i$ and $a^k-b\notin I_i$, one has
$n\in I_j$ for every $j\ne i$.
\end{enumerate}
\end{lemma}
\begin{proof}
This is the arbitrary-product criterion of \cite{TS}. To recall the
argument, test a single coordinate using zeros elsewhere to obtain
(i). If the linear factor is nonzero at $i$, place $(1,1)$ at another
coordinate $j$; its cofactor is $n$, which proves (ii). Conversely,
if one coordinate has a nonzero linear factor, (ii) makes $n$ zero
in every other quotient. At those coordinates either the individual
mpdf condition gives zero cofactor, or the linear factor is zero and
the cofactor is $na^{k(n-1)}=0$. Thus the cofactor vanishes in every
quotient.
\end{proof}

\begin{corollary} \label{cor:dedekind-tame}
Let $D$ be a Dedekind domain, let $0\ne I\subsetneq D$ have prime-ideal
factorization $I=\prod_{i=1}^r P_i^{e_i}$ with distinct maximal ideals
$P_i$, and let $n\ge2$, $k\ge1$. Suppose $nD+I=D$.
Then $I$ is $(kn,n)$-mpdf if and only if one of the following holds:
\begin{enumerate}
\item $r=1$, and with $h=\lceil e_1/n\rceil$ one has
      $h\le k$ and $e_1\le(n-1)h$;
\item $r\ge2$, every $e_i=1$, and $\mu_n(D/P_i)=\{1\}$ for every $i$.
\end{enumerate}
For the ring of integers of a number field, put $q_i=|D/P_i|$.
The root condition in (ii) is equivalent to $\gcd(n,q_i-1)=1$.
\end{corollary}
\begin{proof}
The Chinese remainder theorem and localization give
\[
 D/I\cong\prod_{i=1}^r D/P_i^{e_i}
       \cong\prod_{i=1}^r D_{P_i}/P_i^{e_i}D_{P_i}.
\]
Each $D_{P_i}$ is a discrete valuation ring, and $n$ is a unit in it.
Thus a single factor is governed by Theorem~\ref{thm:tame-dvr}.
In a truncated discrete valuation ring of length $e\ge2$, the choice
$a=0$, $b=\pi^{\lceil e/n\rceil}$ gives equality of the required powers
with a nonzero linear factor. For $e=1$, the quotient is a field $K$,
and the linear factor is forced to vanish precisely when
$\mu_n(K)=\{1\}$: divide by $a^k$ when $a\ne0$, and use $a=1$ to
realize every $n$th root of unity.
If $r\ge2$, Lemma~\ref{lem:product} forces this linear alternative in
every component: the condition $nD+I=D$ makes $n$ a unit modulo every $P_j^{e_j}$,
so $n\notin I_j$ for each component.
This is exactly (ii). The finite-field reformulation follows from the
cyclicity of the multiplicative group.
\end{proof}

For a prime exponent, the tame and ramified local criteria combine
without imposing a condition on the number or norms of the prime ideals.

\begin{theorem}
\label{thm:dedekind-prime}
Let $D$ be a Dedekind domain of characteristic zero, let $p$ be a prime
number, and let $0\ne I\subsetneq D$ have factorization
$I=\prod_{i=1}^rP_i^{e_i}$ with the $P_i$ distinct. Fix $k\ge1$ and put
$E_i=v_{P_i}(p)$ and $h_i=\lceil e_i/p\rceil$.
Then $I$ is $(kp,p)$-mpdf if and only if both of the following conditions
hold:
\begin{enumerate}
\item For each $i$, if $E_i=0$, then $h_i\le k$ and
      $e_i\le(p-1)h_i$; if $E_i>0$, then $e_i\le E_i$ and
      $e_i\le(p-1)h_i$.
\item For every index $i$ such that $e_i\ge2$ or
      $\mu_p(D/P_i)\ne\{1\}$, one has $e_j\le E_j$ for every $j\ne i$.
\end{enumerate}
For $D$ the ring of integers of a number field, the root condition is
$\mu_p(D/P_i)\ne\{1\}$ if and only if $p\mid(|D/P_i|-1)$.
\end{theorem}
\begin{proof}
By the Chinese remainder theorem, $D/I\cong\prod_iV_i/\pi_i^{e_i}V_i$, where
$V_i=D_{P_i}$ is a discrete valuation ring with uniformizer $\pi_i$.
We apply Lemma~\ref{lem:product} with $n=p$.

(a) The $i$th factor satisfies $\MP{k}{p}$ if and only if the inequalities in (i) hold.
For $E_i=0$ this is Theorem~\ref{thm:tame-dvr}. For $E_i>0$ the ring $V_i$ has
characteristic zero and residue characteristic $p$, and this is
Theorem~\ref{thm:ramified-prime}.

(b) The $i$th factor contains $a,b$ with $a^{kp}=b^p$ and $a^k\ne b$ if and only if
$e_i\ge2$ or $\mu_p(D/P_i)\ne\{1\}$.
If $e_i\ge2$, take $a=0$ and $b=\pi_i^{e_i-1}$: since $p(e_i-1)\ge e_i$ we get
$b^p=0$, while $b\ne0$.
If $e_i=1$, the factor is a field $K$. For $a\ne0$ we have $b/a^k\in\mu_p(K)$, and
$a=1$ realizes every element of $\mu_p(K)$; for $a=0$ we get $b=0$.
When $E_i>0$ the residue characteristic is $p$, so $\mu_p(D/P_i)=\{1\}$ and
the condition reduces to $e_i\ge2$.

(c) For $j\ne i$, we have $p=0$ in $V_j/\pi_j^{e_j}V_j$ if and only if $e_j\le E_j$.

Lemma~\ref{lem:product} together with (a), (b), (c) gives (i) and (ii).
For a number field, $(D/P_i)^\times$ is cyclic of order $|D/P_i|-1$, so
$\mu_p(D/P_i)\ne\{1\}$ if and only if $p\mid|D/P_i|-1$.
\end{proof}
\section{All ideals of the integers}\label{sec:integers}

The local valuation criteria specialize to prime-power ideals of
$\mathbb Z$. The global classification then combines these tests with
the compatibility of factor choices across the Chinese remainder
decomposition. The zero ideal is treated separately at the end of the
classification theorem.

\begin{corollary}\label{thm:integer-pp}
Let $p$ be prime, $n\ge2$, and $k,e\ge1$. Put $t=v_p(n)$.
Then $p^e\mathbb Z$ is $(kn,n)$-mpdf precisely in the following cases:
\begin{enumerate}
\item If $t=0$, put $h=\lceil e/n\rceil$. Then $h\le k$ and
$e\le(n-1)h$.
\item If $p$ is odd and $t>0$, then $e\le t$.
\item If $p=2$ and $t>0$, then $e\le t$, or
$t=1$, $e=2$, and $n\ge6$.
\end{enumerate}
\end{corollary}
\begin{proof}
Apply Theorem~\ref{thm:unramified-dvr} to $V=\mathbb Z_{(p)}$, whose
residue field is $\mathbb F_p$, using
$\mathbb Z/p^e\mathbb Z\cong V/p^eV$ and the quotient criterion.
\end{proof}

To apply the product criterion, it remains to identify the local
quotients in which every power equality forces the linear factor to vanish.

\begin{lemma}\label{lem:integer-linear}
Let $p$ be prime, $e,k\ge1$, and $n\ge2$. In $\mathbb Z/p^e\mathbb Z$,
the implication
\[
 a^{kn}=b^n\quad\Longrightarrow\quad a^k=b
\]
holds for every $a,b$ if and only if $e=1$ and $\gcd(n,p-1)=1$.
\end{lemma}
\begin{proof}
For $e\ge2$, take $a=0$ and $b=p^{\lceil e/n\rceil}$.
Then $b^n=0$ modulo $p^e$, while $b\ne0$, since
$\lceil e/n\rceil<e$.
For $e=1$, the multiplicative group of the field $\mathbb F_p$ is cyclic
of order $p-1$. If $\gcd(n,p-1)=1$, raising to the $n$th power is injective,
so the implication holds. If the gcd is larger than one, choose a
nontrivial $n$th root of unity $\zeta\in\mathbb F_p$, and take
$a=1$, $b=\zeta$.
\end{proof}

The local tests and their compatibility now give a finite arithmetic
criterion for every nonzero proper ideal of $\mathbb Z$.

\begin{theorem}\label{thm:integer-all}
Let $N\ge2$ have prime factorization $N=\prod_{p\mid N}p^{e_p}$,
and let $n\ge2$, $k\ge1$.
Then $N\mathbb Z$ is $(kn,n)$-mpdf if and only if both of the following
hold:
\begin{enumerate}
\item Every pair $(p,e_p)$ satisfies the corresponding local condition
      in Corollary~\ref{thm:integer-pp}.
\item For every $p\mid N$ such that
      $e_p\ge2$ or $\gcd(n,p-1)>1$, one has
      $N/p^{e_p}\mid n$.
\end{enumerate}
The zero ideal of $\mathbb Z$ is $(kn,n)$-mpdf for all $n,k$.
\end{theorem}
\begin{proof}
By the Chinese remainder theorem,
$\mathbb Z/N\mathbb Z\cong\prod_{p\mid N}\mathbb Z/p^{e_p}\mathbb Z$.
Apply Lemmas~\ref{lem:product} and \ref{lem:integer-linear}.
The condition that $n=0$ in every other factor is exactly
$\prod_{q\mid N,\,q\ne p}q^{e_q}\mid n$, which is condition (ii).
Finally $(0)$ is prime in $\mathbb Z$.
\end{proof}

\begin{corollary}
Let $n\ge2$, $k\ge1$, and $N\ge2$, and suppose $\gcd(N,n)=1$. If $N$ has at least two distinct
prime divisors, then $N\mathbb Z$ is $(kn,n)$-mpdf if and only if $N$
is squarefree and $\gcd(n,p-1)=1$ for every prime $p\mid N$.
For $N=p^e$, the condition is exactly the tame DVR criterion.
\end{corollary}
\begin{proof}
If $N$ has at least two prime divisors, the integer $N/p^{e_p}>1$ is
coprime to $n$, so it cannot divide $n$. Therefore every component
must satisfy the forced-linear condition of Lemma~\ref{lem:integer-linear}.
The converse follows from the injectivity of the $n$th-power map in
every field component.
\end{proof}

The cubic case illustrates both the extra valuation level allowed by mixed
powers and the restriction on accompanying prime factors.
The squarefree subcase is already \cite[Theorem~5]{YSO}.
Their cdf-primary classification \cite[Theorem~6]{YSO} allows
unrestricted prime powers; the strict condition here imposes the
following bounds, which also depend on $k$.

\begin{corollary}\label{cor:integer-cubic}
Let $N\ge2$ and $k\ge1$. Put $E_1=\{1,2\}$ and
$E_k=\{1,2,4\}$ for $k\ge2$. Then $N\mathbb Z$ is $(3k,3)$-mpdf
if and only if at least one of the following holds:
\begin{enumerate}
\item $N$ is squarefree and no prime divisor of $N$ is congruent to
$1$ modulo $3$;
\item $N=q^e$ or $N=3q^e$, where $q\ne3$ is prime and $e\in E_k$.
\end{enumerate}
\end{corollary}
\begin{proof}
Here $n=3$. At the prime $3$, Corollary~\ref{thm:integer-pp} allows only the exponent $1$.
At a prime $q\ne3$, Theorem~\ref{thm:tame-dvr} gives $h=\lceil e/3\rceil\le k$ and $e\le2h$:
$h=1$ gives $e\in\{1,2\}$, $h=2$ gives $e=4$ (and requires $k\ge2$), and $e\ge5$ is impossible.
So the local condition at $q$ is exactly $e\in E_k$.
Call a prime divisor $q$ of $N$ \emph{special} if $e_q\ge2$ or $\gcd(3,q-1)>1$, that is, $q\equiv1\pmod3$.
The prime $3$ is never special.
Suppose first that $N$ has a special prime $q$. Theorem~\ref{thm:integer-all}(ii) gives
$N/q^{e_q}\mid3$, so $N=q^e$ or $N=3q^e$ with $q\ne3$, and the local condition at $q$
gives $e\in E_k$. This is (ii).
Suppose next that $N$ has no special prime. Then $N$ is squarefree and no prime divisor
of $N$ is congruent to $1$ modulo $3$, which is (i).
Conversely, in case (i) every exponent is $1$, so the local conditions hold
($h=1\le k$ and $1\le 2$ at primes $q\ne3$, and $1\le v_3(3)$ at $3$), and
condition (ii) of Theorem~\ref{thm:integer-all} is vacuous because no prime is special.
In case (ii) the local conditions hold by the choice of $e$, and every special prime
$q$ satisfies $N/q^e\in\{1,3\}$, which divides $3$.
\end{proof}

\begin{example}\label{ex:cubic-comparison}
The ideals $16\mathbb Z$ and $48\mathbb Z$ are $(3k,3)$-mpdf exactly
when $k\ge2$, whereas $12\mathbb Z$ is $(3k,3)$-mpdf for every $k\ge1$.
The ideal $20\mathbb Z$ is never $(3k,3)$-mpdf: the nonreduced component
modulo $4$ can select a nonzero linear factor, and the other component
modulo $5$ does not annihilate $3$. By contrast, every squarefree product
of primes congruent to $2$ modulo $3$, optionally multiplied by $3$,
is $(3k,3)$-mpdf for every $k$.
\end{example}

\begin{remark}\label{rem:scope}
The radical hypotheses in Theorem~\ref{thm:series} cannot simply be
removed. As observed in the polynomial-ascent counterexample of \cite{TS}, the ring
$\mathbb Z/4\mathbb Z$ satisfies $\MP{1}{6}$, but its polynomial ring
does not. More generally the failure holds for every $k$: take $a=1$
and $b=1+2X$. Then $a^{6k}=b^6=1$, whereas
$a^k-b=-2X\ne0$ and $F_{k,6}(a,b)=2+2X\ne0$.
Indeed, $\mathbb Z/4\mathbb Z[X]$ is a square-zero extension of the
domain $\mathbb F_2[X]$ whose kernel $2\,\mathbb Z/4\mathbb Z[X]$ is infinite,
so it violates condition (iii) of Theorem~\ref{thm:squarezero-domain}. Similarly, the ramified DVR theorem treats
prime exponents; it does not assert the same formula for arbitrary
composite exponents in an arbitrarily ramified DVR.
\end{remark}

\end{document}